\documentclass[final,onefignum,onetabnum]{siamart220329}

\usepackage{braket,amsfonts}

\usepackage{array}

\usepackage[caption=false]{subfig}
\usepackage{pgfplots}

\newsiamthm{claim}{Claim}
\newsiamremark{remark}{Remark}
\newsiamremark{hypothesis}{Hypothesis}
\crefname{hypothesis}{Hypothesis}{Hypotheses}

\usepackage{algorithmic}

\usepackage{graphicx,epstopdf}

\Crefname{ALC@unique}{Line}{Lines}

\usepackage{amsopn}

\usepackage{xspace}
\usepackage{bold-extra}
\usepackage[most]{tcolorbox}
\usepackage{showkeys}

\colorlet{texcscolor}{blue!50!black}
\colorlet{texemcolor}{red!70!black}
\colorlet{texpreamble}{red!70!black}
\colorlet{codebackground}{black!25!white!25}

\lstdefinestyle{siamlatex}{%
  style=tcblatex,
  texcsstyle=*\color{texcscolor},
  texcsstyle=[2]\color{texemcolor},
  keywordstyle=[2]\color{texemcolor},
  moretexcs={cref,Cref,maketitle,mathcal,text,headers,email,url},
}

\tcbset{%
  colframe=black!75!white!75,
  coltitle=white,
  colback=codebackground, 
  colbacklower=white, 
  fonttitle=\bfseries,
  arc=0pt,outer arc=0pt,
  top=1pt,bottom=1pt,left=1mm,right=1mm,middle=1mm,boxsep=1mm,
  leftrule=0.3mm,rightrule=0.3mm,toprule=0.3mm,bottomrule=0.3mm,
  listing options={style=siamlatex}
}

\newtcblisting[use counter=example]{example}[2][]{%
  title={Example~\thetcbcounter: #2},#1}

\newtcbinputlisting[use counter=example]{\examplefile}[3][]{%
  title={Example~\thetcbcounter: #2},listing file={#3},#1}

\DeclareTotalTCBox{\code}{ v O{} }
{ 
  fontupper=\ttfamily\color{black},
  nobeforeafter,
  tcbox raise base,
  colback=codebackground,colframe=white,
  top=0pt,bottom=0pt,left=0mm,right=0mm,
  leftrule=0pt,rightrule=0pt,toprule=0mm,bottomrule=0mm,
  boxsep=0.5mm,
  #2}{#1}

\patchcmd\newpage{\vfil}{}{}{}
\begin{tcbverbatimwrite}{tmp_\jobname_header.tex}
\title{A numerical method for solving the generalized tangent vector of hyperbolic systems \thanks{\funding{The first author is funded by the Deutsche Forschungsgemeinschaft (DFG, German Research Foundation) - SPP 2410 Hyperbolic Balance Laws in Fluid Mechanics: Complexity, Scales, Randomness (CoScaRa) within the Project(s) HE5386/26-1 (Numerische Verfahren für gekoppelte Mehrskalenprobleme,525842915),  
(Zufällige kompressible Euler Gleichungen: Numerik und ihre Analysis, 525853336) HE5386/27-1, 
and under HE5386/25-1 Differenzierbare Programmierung für Strömungen mit Diskontinuitäten (513718742).
Furthermore, support received funding from the European Union's Horizon Europe research and innovation programme under the Marie Sklodowska-Curie Doctoral Network Datahyking (Grant No.101072546) is acknowledged.
The second author is funded by Alexander von Humboldt Foundation (Humboldt Research Fellowship Programme for Postdocs).
}}}

\author{Michael Herty \thanks{Chair in Numerical Analysis, IGPM, RWTH Aachen University, Templergraben, 55, D-52062 Aachen, Germany (herty@igpm.rwth-aachen.de)}
\and Yizhou Zhou \thanks{Corresponding author, IGPM, RWTH Aachen University, Templergraben, 55, D-52062 Aachen, Germany (zhou@igpm.rwth-aachen.de)}}

\headers{A numerical method for solving the generalized tangent vector}{Michael Herty, and Yizhou Zhou}
\end{tcbverbatimwrite}
\input{tmp_\jobname_header.tex}

\ifpdf
\hypersetup{ pdftitle={A numerical method for solving the generalized tangent vector of hyperbolic systems} }
\fi

\begin{document}
\maketitle

\begin{tcbverbatimwrite}{tmp_\jobname_abstract.tex}
\begin{abstract}
This work is concerned with the computation of the first-order variation for one-dimensional hyperbolic partial differential equations. In the case of shock waves the main challenge is addressed by developing a numerical method to compute the evolution of the generalized tangent vector introduced by Bressan and Marson (1995) \cite{BM}.
Our basic strategy is to combine the conservative numerical schemes and a novel expression of the interface conditions for the tangent vectors along the discontinuity. Based on this, we propose a simple numerical method to compute the tangent vectors for general hyperbolic systems.
Numerical results are presented for Burgers' equation and a $2\times 2$ hyperbolic system with two genuinely nonlinear fields.

\end{abstract}


\end{tcbverbatimwrite}
\input{tmp_\jobname_abstract.tex}

\section{Introduction}
\label{sec:intro}
This paper deals with the computation of the first-order variation for the system of conservation laws in one dimension:
\begin{align}
&u_t + f(u)_x = 0,\qquad x\in \mathbb{R},\label{original} \\[1mm]
&u|_{t=0} = \bar{u}(x). \label{initialdata}
\end{align}
Here $u\in \mathbb{R}^n$ is the unknown and $f:\mathbb{R}^n \rightarrow \mathbb{R}^n$ is a smooth function of $u$.
We assume that the system \eqref{original} is strictly hyperbolic, and that each characteristic
field is either linearly degenerate or genuinely nonlinear in the sense of Lax \cite{Lax}.

It has been shown that the evolution operator $\mathcal{S}_t: \bar{u}(\cdot) \rightarrow u(t,\cdot) = \mathcal{S}_t \bar{u}(\cdot)$ generated by the conservation law \eqref{original} 
is generically non-differentiable in $L^1(\mathbb{R})$, see e.g. \cite{BM,HHNSS} 
for examples. Therefore, the concept of generalized tangent vector (GTV) has been proposed~\cite{BM} for piecewise Lipschitz continuous solution $u(t,\cdot)$ with a finite number of $N$  discontinuities. 
The differential is described by measures with an absolutely continuous part $v(t,\cdot) \in L^1(\mathbb{R})$ and a singular part. The latter is determined as a variation of the shift
$\xi_{\alpha}(t)$ at each jump location $x_{\alpha}(t)$ ($\alpha=1,2,...,N$). It has been shown (\cite{BM}) that $\xi_{\alpha}$ satisfies an ordinary differential equation and $v$ solves the linear equation 
\begin{align}\label{intro-eqv}
v_t + A(u)v_x + [DA(u)v]u_x = 0,
\end{align}
outside the points $x_{\alpha}(t)$ where $u(t,\cdot)$ suffers a discontinuity. Along each line $x = x_{\alpha}(t)$, $v$ satisfies the interface conditions for $t>0$
\begin{align}
&\Psi_{\alpha}(v^-,v^+;u^-,u^+,\xi_{\alpha})=0,\qquad \alpha=1,2,...,N.
\label{intro-eqv3}
\end{align}
Here, $v^+=\lim_{x\rightarrow x_\alpha^+} v(t,x)$ and $v^-=\lim_{x\rightarrow x_\alpha^-} v(t,x)$ represent the left and right limits of $v$ at the  point of discontinuity $x=x_\alpha(t)$. Note that the interface conditions depend on $\xi_{\alpha}$ as well as the limits $u^-$ and $u^+$. We refer to Section \ref{sec:2} for the explicit formula \eqref{theorem2.1:eq3}.

At this point, we briefly review other theoretical results treating differentiability and related questions of hyperbolic conservation laws. 
In \cite{MR1816648}, the Lipschitz continuous dependence of $\mathcal{S}_t$ in $L^1$ on the initial data has been established. Furthermore, an adjoint calculus and optimality conditions are derived in \cite{MR2311524}.
In the scalar 1D case, results on differentiability have been extended from directional  to Fréchet-type differentiability \cite{ulbrich2001optimal,MR2397976,MR1939870,MR2020647,MR4241509}  in the case of piecewise $C^1$ initial data. Moreover, a sensitivity and adjoint calculus is developed for both initial and initial-boundary value problems. An alternative variational calculus, referred to as shift-differentiability  based on horizontal shifts of the graph, is developed in \cite{MR1422538} for the scalar case and in \cite{MR1739381} for the system's case. We refer to the references therein for further details. 

The main goal of the present work is to propose a (simple) numerical  method for the computation of the generalized tangent vector $(v,\xi_{\alpha})$ introduced in \cite{BM}. In particular, our main concern is to develop a strategy dealing with interface conditions \eqref{intro-eqv3} for the GTV. 
Our work is motivated by numerical methods for hyperbolic conservation laws. Therein, computing  discontinuous solutions of \eqref{original}, can be obtained using  conservative schemes that  preserve the Rankine–Hugoniot condition
\begin{equation*}
s[u] = [f(u)]\qquad \text{at the jump position} ~~x=x_{\alpha}(t).
\end{equation*}
Here $s=x_{\alpha}'(t)$ is the shock speed and  $[u]=u^+-u_-$.
This strategy can also be used to compute the generalized tangent vector. The interface conditions \eqref{intro-eqv3} 
can be regarded as the formal differentiation of the Rankine–Hugoniot condition in the case of a single discontinuity of a single family, see \cite{MR3881242,HHNSS}.
Namely, 
$$
\dot s[u] + s[\dot u] = [\dot f] = [A(u)\dot u].
$$
Here the dot means the formal differentiation with respect to the initial data.
Motivated by Lemma 2.1 in \cite{BM}, we write $\dot s=\xi_{\alpha}'$ and $\dot u = w = v+\xi_{\alpha} u_x$ to conclude
\begin{equation}\label{BC-introduction} 
\xi_{\alpha}'[u] + s[w] = [A(u)w]. 
\end{equation}
In the present paper, we rigorously verify the above formal derivation in the framework of \cite{BM} and present a method for their numerical computation. 
To illustrate the method, we first deal with a simple case of only a single jump discontinuity in the original problem \eqref{original}-\eqref{initialdata}.
In this case, we use the formula \eqref{BC-introduction} within  conservative schemes. The interface conditions \eqref{intro-eqv3} are shown to be preserved without additional numerical treatment. 
For the case of multiple discontinuities, we assume, that there is no interaction between two discontinuities. Finally, we present a simple way to combine the results for different jumps and obtain the generalized tangent vector on the whole domain. 

Numerical methods for discretizing the sensitivity equation in hyperbolic conservation laws have been an active research topic. These methods are particularly relevant to numerical computations in optimal control problems for hyperbolic conservation laws. In \cite{MR3881242,MR4102600,MR2645781,GUINOT200961}, the authors develop a numerical method to compute the sensitivity equation by adding delta-type source terms. We remark that our method is similar to this series of work in computing Riemann problems. However, our method is based on the rigorous theoretical result of the GTV in \cite{BM}. Unlike these works, our numerical strategy for the GTV does not need piecewise constant assumption and can be extended to any conservative schemes. More detailed comparisons can be found in Remark \ref{Remark4.2}.
In \cite{HHNSS}, the authors develop an algorithmic differentiation framework for computing the sensitivity equation. 
In the works \cite{GU1,GU2}, the convergence of a modified Lax-Friedrichs scheme for the sensitivity equation has been demonstrated even in the presence of shocks. Additionally, convergence results for state schemes that satisfy a discrete one-sided Lipschitz condition, along with their corresponding adjoint schemes, are presented in \cite{ulbrich2001optimal,SchaeferAguilarSchmittUlbrichMoos2019,SchaeferAguilarUlbrich2021}, while implicit-explicit methods are discussed in \cite{MR2891921,HP}.
In \cite{HP}, a relaxation approximation is exploited in the tangent vector calculus, and corresponding schemes are derived.  Further examples include finite volume, Lagrangian methods, and vanishing viscosity approach in \cite{MR3828251,MR3238137,MR3274869,MR3043133}. We would like to stress the fact that the sensitivity analysis of hyperbolic problems has a rich field of application in engineering \cite{Giles1,LiardZuazua,MR2095284,MR3101096,MR1803927}.

The paper is organized as follows. In Section \ref{sec:2}, we review the theoretical results of the generalized tangent vector. Section \ref{sec:3} is devoted to deduction of interface conditions and two illustrative examples are presented. In Section \ref{sec:4}, we give a numerical method for computing the generalized tangent vector, and more numerical experiments are shown in Section \ref{sec:5}. In Appendix \ref{AppendixA}, we give some details of the proof in Section \ref{sec:3}.

\section{Preliminaries}
\label{sec:2}
This section briefly recalls the theoretical results of the generalized tangent vector (GTV) developed in \cite{BM}.

First, we introduce some basic notations. For the flux term in \eqref{original}, denote the Jacobian matrix $A(u):=f_u(u)$. By the strictly hyperbolic assumption, it has $n$ real distinct eigenvalues $\lambda_1(u) < \lambda_2(u) < \cdots < \lambda_n(u)$. For any two vectors $u,v\in \mathbb{R}^n$, we define the matrix 
\begin{equation}\label{barA}
\bar{A}(u,v) = \int_0^1 A(\theta u + (1-\theta)v) d\theta.
\end{equation}
For $i = 1,...,n$, the $i$-th eigenvalue, left eigenvector and right eigenvector of $\bar{A}(u, v)$ will be denoted by $\lambda_i(u, v)$, $l_i(u, v)$ and $r_i(u, v)$ respectively. We assume that the ranges of the eigenvalues $\lambda_i(u, v)$ do not overlap for any $(u, v)$. Because of the regularity of $A$ and the strict hyperbolic assumption, $\lambda_i(u, v)$, $l_i(u, v)$ and $r_i(u, v)$ are all $C^1$ functions.

As in \cite{BM}, we only consider the piecewise Lipschitz continuous solution $u(t,x)$ with $N$ discontinuities. For each $\alpha = 1, ..., N$, let $x_{\alpha}(t)$ be the position of the $\alpha$-th discontinuity of $u(t,x)$ at time $t$. Assume that the discontinuity at $x=x_{\alpha}(t)$ corresponds to the $k_{\alpha}$-th characteristic
family, the Rankine-Hugoniot and the entropy admissibility conditions imply
\begin{align*}
&x_{\alpha}'(t)=\lambda_{k_{\alpha}}(u^+,u^-),
\\[1mm]
&\langle l_i(u^+,u^-), u^+-u^- \rangle = 0,\qquad \forall ~i\neq k_{\alpha},\\[1mm]
&\lambda_{k_{\alpha}}(u^+,u^+)\leq \lambda_{k_{\alpha}}(u^+,u^-) \leq 
\lambda_{k_{\alpha}}(u^-,u^-).
\end{align*}
Here $u^+=\lim_{x\rightarrow x_\alpha^+} u(t,x)$ and $u^-=\lim_{x\rightarrow x_\alpha^-} u(t,x)$ represent the left and right limits of $u$ at the discontinuity  $x=x_{\alpha}(t)$.

In studying the first-order variation of \eqref{original}, some regularity conditions for the perturbed solution $u^{\epsilon}$ should be given. Namely,
\begin{definition}
Let $u=u(x)$ be a piecewise Lipschitz continuous function
on $\mathbb{R}$ with $N$ simple discontinuities. A path $\gamma$ is a Regular Variation (R.V.) for $u$ if, for $\epsilon \in [0, \epsilon_0]$, all functions $u^\epsilon=\gamma(\epsilon)$ are piecewise Lipschitz continuous, with jumps at points $x_1^{\epsilon} < ... < x_N^{\epsilon}$ depending continuously on $\epsilon$. Moreover, the Lipschitz constant of $u$ is independent of $\epsilon$.
\end{definition}

With these preparations, the main result in \cite{BM} can be stated: 
\begin{theorem}
Let $u=u(t,x)$ be a piecewise
Lipschitz continuous solution of \eqref{original} with $N$ simple discontinuities. 
Let $(\bar{v},\bar{\xi})\in L^1(\mathbb{R})\times \mathbb{R}^N$ be a tangent vector to the initial data $u(0,x)=\bar{u}(x)$, generated by the R.V. $\epsilon \mapsto \bar{u}^{\epsilon}$, and call $u^{\epsilon} = u^{\epsilon}(t,x)$ the solution of \eqref{original} with initial condition $\bar{u}^{\epsilon}$.
Then there exists $\tau_0 > 0$ such that, for all $t \in [0,\tau_0]$, the path $\epsilon \mapsto u^{\epsilon}(t,\cdot)$
is a R.V. for $u(t,\cdot)$, generating the tangent vector $(v(t), \xi(t))$. The vector $(v(t), \xi(t))$ is the unique solution of the problem:
\begin{align}
&v_t + A(u)v_x + [DA(u)v]u_x = 0,\qquad v(0,x)=\bar{v}(x)\label{theorem2.1:eq1} 
\end{align}
outside the discontinuities of $u$, while, for $\alpha = 1,...,N$, 
\begin{align}
&\xi'_{\alpha} = D\lambda_{k_\alpha}(u^+,u^-)\cdot (v^++\xi_{\alpha}u_x^+,v^-+\xi_{\alpha}u_x^-),\qquad \xi_{\alpha}(0)=\bar{\xi}_\alpha, \label{theorem2.1:eq2} \\[2mm]
&\big \langle Dl_{i}(u^+,u^-)\cdot (v^++\xi_{\alpha}u_x^+,v^-+\xi_{\alpha}u_x^-), ~u^+-u^- \big \rangle \nonumber \\[1mm]
&\qquad + \big \langle l_{i}(u^+,u^-),  (v^++\xi_{\alpha}u_x^+ - v^- -\xi_{\alpha}u_x^- \big \rangle = 0,\qquad \forall i \neq k_{\alpha}
\label{theorem2.1:eq3}
\end{align}
along each line $x = x_{\alpha}(t)$ where $u(t,\cdot)$ suffers a discontinuity, in the $k_{\alpha}$-th characteristic direction.
\end{theorem}

This theorem provides the evolution equations for $v=v(t,x)$ and $\xi=\xi(t)$. Note that the formula \eqref{theorem2.1:eq3} is just the aforementioned interface condition \eqref{intro-eqv3}.
The first-order variation $u_\epsilon$ can be computed from 
\begin{equation}\label{construction}
u_\epsilon = u+\epsilon v + \sum_{\xi_{\alpha}<0} (u(x_{\alpha}^+)-u(x_{\alpha}^-))\chi_{[x_{\alpha} + \epsilon \xi_{\alpha}, x_{\alpha}]} - \sum_{\xi_{\alpha}>0} (u(x_{\alpha}^+)-u(x_{\alpha}^-))\chi_{[x_{\alpha}, x_{\alpha} + \epsilon \xi_{\alpha}]}.
\end{equation}
Here $u(x_{\alpha}^\pm)=\lim_{x\rightarrow x_\alpha^\pm} u(t,x)$ and $\chi_{[a,b]}=\chi_{[a,b]}(x)$ is the characteristic function of the interval $[a,b]$. The construction of $u^{\epsilon}$ can be understood in the following way: starting with $u$, adding $\epsilon v$ and then shifting by $\epsilon \xi_{\alpha}$ each point $x_{\alpha}$ where $u$ has a jump.

\section{Interface conditions for the GTV}
\label{sec:3}

We aim to present a numerical method for solving \eqref{theorem2.1:eq1}-\eqref{theorem2.1:eq3}. Recall that the equation \eqref{theorem2.1:eq1} is only valid outside the shock position $x=x_{\alpha}(t)$, i.e., $(t,x)\not=(t,x_\alpha(t))$. Along the discontinuity $x=x_{\alpha}(t)$, the ordinary differential equation \eqref{theorem2.1:eq2} 
\begin{align*}
&\xi'_{\alpha} = D\lambda_{k_\alpha}(u^+,u^-)\cdot (v^++\xi_{\alpha}u_x^+,v^-+\xi_{\alpha}u_x^-)
\end{align*}
and the interface conditions \eqref{theorem2.1:eq3} 
\begin{align*}
&\big \langle Dl_{i}(u^+,u^-)\cdot (v^++\xi_{\alpha}u_x^+,v^-+\xi_{\alpha}u_x^-), ~u^+-u^- \big \rangle \\[1mm]
&\qquad + \big \langle l_{i}(u^+,u^-),  (v^++\xi_{\alpha}u_x^+ - v^- -\xi_{\alpha}u_x^- \big \rangle = 0,\qquad \forall i \neq k_{\alpha}
\end{align*}
hold true. Next, we will present new expressions for  \eqref{theorem2.1:eq2} and \eqref{theorem2.1:eq3}, that allows for a simple numerical treatment in Section \ref{sec:4}.

Recall the proof in \cite{BM} where the variables 
\begin{equation}\label{w-variable}
w_{\alpha}=v+\xi_{\alpha} u_x,\qquad y=x-x_{\alpha}(t)
\end{equation}
are introduced to analyze the ${\alpha}$-th jump in a neighborhood of $x_{\alpha}(t)$. 
The equation for $w_{\alpha}=w_{\alpha}(y,t)$ reads
\begin{align*}
    &\partial_t w_{\alpha} + [A(u) - x_{\alpha}'I] \partial_y w_{\alpha} + DA(u) (\partial_y u) w_{\alpha} - \xi_{\alpha}' \partial_y u = 0.
\end{align*}
Moreover, the relations \eqref{theorem2.1:eq2} and \eqref{theorem2.1:eq3} become
\begin{align}
    &\xi_{\alpha}'(t)=D \lambda_{k_{\alpha}}(u^+,u^-)(w_{\alpha}^{+},w_{\alpha}^{-}) \label{BCxi}\\[2mm]
    &\langle l_i(u^{+},u^{-}), w_{\alpha}^{+}-w_{\alpha}^{-} \rangle + \langle D l_i(u^{+},u^{-})(w_{\alpha}^{+},w_{\alpha}^{-}), u^{+}-u^{-} \rangle = 0,\qquad i\neq k_{\alpha}. \label{BCw}
\end{align}
Assume that the solution is sufficiently smooth except for the position $y=0$. We can rewrite the equation to the conservative form
\begin{equation}\label{11eq5}
    \partial_t w_{\alpha} + \partial_y [A(u)w_{\alpha} - x_{\alpha}'w_{\alpha} - \xi_{\alpha}' u]=0.
\end{equation}
Based on these preparations, we state
\begin{theorem}\label{theorem3.1}
    The equation \eqref{BCxi} is equivalent to the relation:
    \begin{equation}\label{eq-2.1}
    \xi_{\alpha}' = \left\langle l_{k_{\alpha}}(u^+,u^-), ~[A(u^+) - \bar{A}(u^+,u^-)] w^+_{\alpha} + [\bar{A}(u^+,u^-) - A(u^-)] w^-_{\alpha}\right\rangle.
    \end{equation}
Here $l_{k_{\alpha}}(u^+,u^-)$ is the  left eigenvector of $\bar{A}(u^+,u^-)$ satisfying $\langle l_{k_{\alpha}}(u^+,u^-), u^+-u^- \rangle=1$.
Moreover, the interface condition \eqref{BCw} is satisfied under the following relation:
\begin{equation}\label{lem11eq2}
x_{\alpha}' w^+_{\alpha} - A(u^+)w^+_{\alpha} + \xi_{\alpha}' u^+ = x_{\alpha}' w^-_{\alpha} -A(u^-)w^-_{\alpha} + \xi_{\alpha}' u^-,
\end{equation}    
which implies that the flux in \eqref{11eq5} is continuous at $y=0$.
\end{theorem}

In order to prove this theorem, we need the following Lemma.

\begin{lemma}\label{lemma12}
If $f\in C^2(G)$, the following relations hold for each $i=1,2,...,n$:
\begin{align*}
&\frac{\partial l_i}{\partial u^{+}} (\lambda_{k_{\alpha}} - \lambda_i)(u^+-u^-) + l_i [A(u^+) - \bar{A}(u^+,u^-)]
=\frac{\partial \lambda_i}{\partial u^+} \langle l_i, u^+-u^-\rangle, \\[2mm]
&\frac{\partial l_i}{\partial u^{-}} (\lambda_{k_{\alpha}} - \lambda_i)(u^+-u^-) + l_i [ \bar{A}(u^+,u^-)-A(u^-)]
=\frac{\partial \lambda_i}{\partial u^-} \langle l_i, u^+-u^-\rangle.
\end{align*}
Here $\lambda_i=\lambda_i(u^+,u^-)$ is the $i$-th eigenvalue of $\bar{A}(u^+,u^-)$ and $l_i=l_i(u^+,u^-)$ is the corresponding (left) eigenvector.
\end{lemma}
We leave the proof of this lemma to Appendix \ref{AppendixA}. We now prove Theorem \ref{theorem3.1} as follows.

\begin{proof}
Taking $i=k_\alpha$ in Lemma \ref{lemma12} and using $\langle l_{k_\alpha}, u^+-u_-\rangle = 1$, we have 
$$
\frac{\partial \lambda_{k_{\alpha}}}{\partial u^+}
=
l_{k_{\alpha}} [A(u^+) - \bar{A}(u^+,u^-)],\qquad
\frac{\partial \lambda_{k_{\alpha}}}{\partial u^-}
=
l_{k_{\alpha}} [\bar{A}(u^+,u^-)-A(u^-)].
$$
Then it follows that 
\begin{align*}
\xi_{\alpha}' 
= &~ D \lambda_{k_{\alpha}}(u^+,u^-)(w_{\alpha}^{+},w_{\alpha}^{-}) 
= \Big\langle \frac{\partial \lambda_{k_{\alpha}}}{\partial u^+},~ w_{\alpha}^{+}\Big\rangle + \Big\langle\frac{\partial \lambda_{k_{\alpha}}}{\partial u^-}, ~w_{\alpha}^{-} \Big\rangle \\[2mm]
= &~ \left\langle l_{k_{\alpha}} [A(u^+) - \bar{A}(u^+,u^-)], ~w_{\alpha}^+\right\rangle
+ \left\langle l_{k_{\alpha}}[\bar{A}(u^+,u^-)-A(u^-)], ~w_{\alpha}^-\right\rangle.
\end{align*}
This is the equation \eqref{eq-2.1}. Furthermore, we write \eqref{lem11eq2} as
\begin{align*}
\xi_{\alpha}'(u^+-u^-) = A(u^+)w_{\alpha}^+ - A(u^-)w_{\alpha}^- - x_{\alpha}'(w_{\alpha}^+-w_{\alpha}^-).    
\end{align*}
Multiplying $l_i$ on the left and using the fact $\langle l_i, u^+-u^-\rangle=0$ for any $i\neq k_{\alpha}$, we obtain
\begin{align}
    0 =&~ \langle l_i, A(u^+)w_{\alpha}^+-A(u^-)w_{\alpha}^- - x_{\alpha}'(w_{\alpha}^+-w_{\alpha}^-) \rangle,\nonumber \\[2mm]
     =&~ \langle l_i, [A(u^+) - \bar{A}(u^+,u^-)] w_{\alpha}^+ +
    [\bar{A}(u^+,u^-) - A(u^-)] w_{\alpha}^- \rangle \nonumber \\[2mm]
    & +\langle l_i, \bar{A}(u^+,u^-) 
    (w_{\alpha}^+ - w_{\alpha}^-) - x_{\alpha}'(w_{\alpha}^+-w_{\alpha}^-) \rangle \nonumber. 
\end{align}
Since $l_i \bar{A}(u^+,u^-) = \lambda_i l_i $ and 
$x'_\alpha = \lambda_{k_\alpha}$, we have
\begin{align}
    0=&~ \langle l_i, [A(u^+) - \bar{A}(u^+,u^-)] w_{\alpha}^+ +
    [\bar{A}(u^+,u^-) - A(u^-)] w_{\alpha}^- \rangle \label{proofth31:eq1}\\[1mm]
    &~+ (\lambda_i-\lambda_{k_\alpha}) \langle l_i, 
    w_{\alpha}^+ - w_{\alpha}^-\rangle .\nonumber 
\end{align}
In order to simplify this equation, we take $i\neq k_\alpha$ in Lemma \ref{lemma12} and use the relation $\langle l_i, u^+-u^-\rangle=0$ to get
\begin{align*}
0=&~\frac{\partial l_i}{\partial u^{+}} (\lambda_{k_{\alpha}} - \lambda_i)(u^+-u^-) + l_i [A(u^+) - \bar{A}(u^+,u^-)],\\[1mm]
0=&~\frac{\partial l_i}{\partial u^{-}} (\lambda_{k_{\alpha}} - \lambda_i)(u^+-u^-) + l_i [ \bar{A}(u^+,u^-)-A(u^-)].
\end{align*}
Substituting these into \eqref{proofth31:eq1}, we obtain
\begin{align}
    0=&~ \Big\langle \frac{\partial l_i}{\partial u^{+}} (u^+-u^-), w_{\alpha}^+ \Big\rangle + 
    \Big\langle \frac{\partial l_i}{\partial u^{-}} (u^+-u^-), w_{\alpha}^- \Big\rangle
    + \langle l_i, 
    w_{\alpha}^+ - w_{\alpha}^-\rangle .\nonumber 
\end{align}
This is the interface condition \eqref{BCw} and we complete the proof of Theorem \ref{theorem3.1}.
\end{proof}

Next, we illustrate Theorem \ref{theorem3.1} through two examples. 

\textbf{Example 1 (Burgers equation):}
The Riemann problem for the Burgers equation reads as
$$
u_t + \left(\frac{u^2}{2}\right)_x = 0
$$
with the initial data 
$$
u(x,0) = \left\{\begin{array}{cc}
    u^-, & x<0, \\[2mm]
    u^+, & x\geq0.
\end{array}\right.
$$
Here $u^-$ and $u^+$ are two constants. Since there is only one shock, we drop the subscription $\alpha$ and write the generalized tangent vector in \eqref{BCxi} and \eqref{11eq5} as $(\xi,w)$. By a direct computation, we have 
$$
A(u^+) = u^+,\quad A(u^-) = u^-,\quad \bar{A}(u^+,u^-) = \frac{u^++u^-}{2},\quad l(u^+,u^-) = \frac{1}{u^+-u^-}.
$$ 
Then the equation \eqref{BCxi} for $\xi$ reads
$$
\xi' = D \lambda(u^+,u^-)(w^{+},w^{-}) = \lim_{\epsilon\rightarrow 0}\frac{(u^+ + \epsilon w^+)+(u^- +\epsilon w^-) }{2\epsilon} = \frac{w_-+w_+}{2}.
$$
Now we check that the new expression \eqref{eq-2.1} gives the same relation. Indeed, we compute from \eqref{eq-2.1} that
\begin{align*}
    \xi' = \frac{1}{u^+-u^-} \left(u^+ w_+ - \frac{u^++u^-}{2} w_+ + \frac{u^++u^-}{2} w_- - u^- w_- \right) = \frac{w_-+w_+}{2}.
\end{align*}
In the scalar case, the last condition \eqref{BCw} does not exist. We consider the following example for the computation of both \eqref{BCxi} and \eqref{BCw}. 

\textbf{Example 2:}
In order to illustrate the main issues for the system case, we consider a simple $2\times 2$ hyperbolic system 
\begin{eqnarray}\label{psystem}
    \begin{aligned}
        &\partial_t\rho + \partial_x q = 0,\\[2mm]
        &\partial_t q + \partial_x p(\rho) = 0,
    \end{aligned}
\end{eqnarray}
where $p(\rho)=\kappa \rho^\gamma$. Consider the Riemann problem 
$$
(\rho,q)(x,0) = \left\{\begin{array}{cc}
    (\rho^-,q^-), & x<0, \\[2mm]
    (\rho^+,q^+), & x\geq0
\end{array}\right.
$$
with $u^-=(\rho^-,q^-)$ and $u^+=(\rho^+,q^+)$ two constant states.
We compute coefficient matrices
$$
A(u^-)=
\begin{pmatrix}
    0 & 1\\[1mm]
    p'(\rho^-) & 0
\end{pmatrix},\quad
A(u^+)=
\begin{pmatrix}
    0 & 1\\[1mm]
    p'(\rho^+) & 0
\end{pmatrix},\quad 
\bar{A}(u^+,u^-)=
\begin{pmatrix}
    0 & 1\\[1mm]
    \bar{p}'(\rho^+,\rho^-) & 0
\end{pmatrix}
$$
with $\bar{p}'(\rho^+,\rho^-) = (p(\rho^+)-p(\rho^-))/(\rho^+-\rho^-)$ and eigenvalues 
\begin{eqnarray}
    \lambda_1(u^+,u^-) = - \sqrt{\bar{p}'(\rho^+,\rho^-)},\qquad \lambda_2(u^+,u^-) = \sqrt{\bar{p}'(\rho^+,\rho^-)}.
\end{eqnarray}
The 1-shock is characterized by
\begin{align}\label{1-shock}
    q^+-q^-= - \sqrt{( p(\rho^+) -  p(\rho^-))(\rho^+-\rho^-)},\qquad \rho^+ > \rho^-,
\end{align}
while the 2-shock is characterized by
\begin{align}\label{2-shock}
    q^+-q^-= - \sqrt{( p(\rho^+) -  p(\rho^-))(\rho^+-\rho^-)},\qquad \rho^+<\rho^-.
\end{align}
For simplicity, we show the case where only the 1-shock occurs. Namely, $(\rho_+,q_+)$ and $(\rho_-,q_-)$ are given by \eqref{1-shock}.

Next, we derive the interface conditions \eqref{BCxi}-\eqref{BCw} and compare them to the new expression \eqref{eq-2.1}-\eqref{lem11eq2}. Firstly, we denote $w^+=(w_\rho^+,w_q^+)$ and $w^-=(w_\rho^-,w_q^-)$. A direct computation of \eqref{BCxi} shows that 
\begin{align}
    \xi' =&~ \lim_{\epsilon \rightarrow 0}\frac{1}{\epsilon}  \Big(\lambda_{1}(u^++\epsilon w_{\alpha}^{+},u^- + \epsilon w_{\alpha}^{-}) - \lambda_{1}(u^+,u^-)\Big) \nonumber\\[1mm]
    =&~ - \Big( \frac{\partial \sqrt{\bar{p}'(\rho^+,\rho^-)}}{\partial \rho^+}w_\rho^+ + \frac{\partial \sqrt{\bar{p}'(\rho^+,\rho^-)}}{\partial \rho^-}w_\rho^- \Big) \nonumber \\[1mm]
    =&~ 
      -\frac{1}{2}\frac{[p'(\rho^+) - \bar{p}'(\rho^+,\rho^-)]w_{\rho}^+ + 
       [\bar{p}'(\rho^+,\rho^-) - p'(\rho^-)]w_\rho^-}{(\rho^+-\rho^-)\sqrt{\bar{p}'(\rho^+,\rho^-)}}. \label{example2-eq1}
\end{align}
For the interface condition \eqref{BCw}, we compute 
$l_2(u^+,u^-) = (\sqrt{\bar{p}'(\rho^+,\rho^-)},~1)$
and thereby
\begin{align*}
D l_2(u^{+},u^{-})(w^{+},w^{-})
=&~
\lim_{\epsilon \rightarrow 0}\frac{1}{\epsilon}\left(
\sqrt{\bar{p}'(\rho^++\epsilon w_\rho^+,\rho^-+\epsilon w_\rho^-)}
-\sqrt{\bar{p}'(\rho^+,\rho^-)}, ~~0\right)\\[1mm]
=&~\Big(\frac{\partial \sqrt{\bar{p}'(\rho^+,\rho^-)}}{\partial \rho^+}w_\rho^+ + \frac{\partial \sqrt{\bar{p}'(\rho^+,\rho^-)}}{\partial \rho^-}w_\rho^-,~~0\Big)\\[1mm]
=&~ \Big(\frac{1}{2}\frac{[p'(\rho^+) - \bar{p}'(\rho^+,\rho^-)]w_{\rho}^+ + 
       [\bar{p}'(\rho^+,\rho^-) - p'(\rho^-)]w_\rho^-}{(\rho^+-\rho^-)\sqrt{\bar{p}'(\rho^+,\rho^-)}},~~0\Big).
\end{align*}
Then the interface condition \eqref{BCw} reads as
\begin{align}
0=&~\langle l_2(u^{+},u^{-}), w^{+}-w^{-} \rangle + \langle D l_2(u^{+},u^{-})(w^{+},w^{-}), u^{+}-u^{-} \rangle \nonumber \\[2mm]
=&~ \sqrt{\bar{p}'(\rho^+,\rho^-)}(w^{+}_\rho-w^{-}_\rho) + (w^{+}_q-w^{-}_q) \nonumber\\[1mm]
&~ +\frac{1}{2} \frac{[p'(\rho^+) - \bar{p}'(\rho^+,\rho^-)]w_{\rho}^+ + 
       [\bar{p}'(\rho^+,\rho^-) - p'(\rho^-)]w_\rho^-}{\sqrt{\bar{p}'(\rho^+,\rho^-)}} \nonumber \\
 =&~ (w^{+}_q-w^{-}_q) +\frac{1}{2} \frac{p'(\rho^+)  w_{\rho}^+ - p'(\rho^-)w_\rho^-}{\sqrt{\bar{p}'(\rho^+,\rho^-)}} + \frac{1}{2}\sqrt{\bar{p}'(\rho^+,\rho^-)}(w^{+}_\rho-w^{-}_\rho) .  \label{example2-eq2}  
\end{align}
On the other hand, we use the expression in Theorem \ref{theorem3.1}. 
The left eigenvector $l_1=l_1(u^+,u^-)$ satisfying $\langle l_1, u^+-u^-\rangle = 1$ is given by
$$
l_1(u^+,u^-) 
=  \frac{1}{2(\rho^+-\rho^-)\sqrt{\bar{p}'(\rho^+,\rho^-)}}\Big(-\sqrt{\bar{p}'(\rho^+,\rho^-)},~1\Big).
$$
Thus it follows from \eqref{eq-2.1} that
\begin{align*}
    \xi' =&~ \left\langle l_1(u^+,u^-), ~[A(u^+) - \bar{A}(u^+,u^-)] w^+ + [\bar{A}(u^+,u^-) - A(u^-)] w^-\right\rangle \\[1mm]
    =&~ -
      \frac{1}{2}\frac{[p'(\rho^+) - \bar{p}'(\rho^+,\rho^-)]w_{\rho}^+ + 
       [\bar{p}'(\rho^+,\rho^-) - p'(\rho^-)]w_\rho^-}{(\rho^+-\rho^-)\sqrt{\bar{p}'(\rho^+,\rho^-)}},
\end{align*}
which is just the formula \eqref{example2-eq1}. At last, 
\eqref{lem11eq2} implies
$$
-\sqrt{\bar{p}'(\rho^+,\rho^-)}
\begin{pmatrix}
w^+_{\rho}-w^-_{\rho}\\[1mm]
w^+_{q}-w^-_{q}
\end{pmatrix}
=
\begin{pmatrix}
w^+_{q}-w^-_{q}\\[1mm]
p'(\rho^+)w^+_{\rho}-p'(\rho^-)w^-_{\rho}
\end{pmatrix}
+ \xi'
\begin{pmatrix}
\rho^+-\rho^-\\[1mm]
q^+-q^-
\end{pmatrix}.
$$
Eliminating $\xi'$ from this equation and using \eqref{1-shock} yield \eqref{example2-eq2}.

\section{Numerical method}\label{sec:4}
In the previous section, an explicit expression for the interface conditions at $x=x_\alpha(t)$ has been introduced. We will show that this expression can be used to design the numerical method. 

\subsection{Single shock case (N=1)} \label{sec:4.1}
At first, we illustrate the basic idea through a simple case where $N=1$. Namely, there is only one discontinuous point for the solution $u$ to the original equation \eqref{original}. For simplicity, we drop the subscription $\alpha$ for $\xi_{\alpha}$ and $w_{\alpha}$.
In the equation \eqref{11eq5}, the interface is fixed at $y=0$. Now we use the original spatial variable $x$ and rewrite \eqref{11eq5} as:  
\begin{equation}\label{eq:W}
    \partial_t w + \partial_x [A(u)w - \xi' u]=0,
    \qquad x\in (-\infty,x_{\alpha})\cup(x_{\alpha},+\infty).
\end{equation} 
Denoting $F = A(u)w - \xi'u$, we know that the weak solution satisfies
$$
\int_{D} \Big( \phi_t^T w + \phi^T_x F \Big)~ dx dt = 0,\qquad \forall \phi\in C_0^{\infty}(D)
$$
with the domain $D\in\{x\in\mathbb{R},t>0\}$. Consider $D^-=\{(x,t)\in D,~x<x_\alpha(t)\}$ and $D^+=\{(x,t)\in D,~x>x_\alpha(t)\}$.
We compute
\begin{align*} 
0=&\int_{D^-} (\phi^T w)_t + (\phi^T F)_x~ dx dt + \int_{D^+} (\phi^T w)_t + (\phi^T F)_x~ dx dt\\
=&\int_{\partial D^-} (-\phi^T w)~ dx + (\phi^T F )~ dt + \int_{\partial D^+} (-\phi^T w)~ dx + (\phi^T F )~ dt\\
=& \int_{x=x_\alpha(t)} -\phi^T (w^- - w^+)~ dx + \phi^T (F^- - F^+)~ dt
\end{align*}
with $F^{+} = A(u^{+})w^{+} - \xi'u^{+}$ and $F^{-} = A(u^{-})w^{-} - \xi'u^{-}$.
Since $\phi$ is arbitrary and $dx/dt = x_{\alpha}'$, we conclude that 
$$
x'(w^+-w^-) = A(u^+)w^+ - \xi' u^+ - A(u^-)w^- + \xi' u^-,
$$
which is just \eqref{lem11eq2}. 

According to the classical theory of numerical schemes for conservation laws \cite{Le}, the relation \eqref{lem11eq2}  can be preserved along the discontinuity $x=x_\alpha(t)$ if we use a conservative numerical scheme to solve equation \eqref{eq:W} in the domain $D$. Moreover, Theorem \ref{theorem3.1} implies that the interface condition \eqref{BCw} is satisfied. 
Motivated by this argument, a numerical treatment is  straightforward. 

Before we discuss our algorithm, we introduce some preparations.
Consider a mesh on $\mathbb{R}$ consisting of cells $I_i = (x_{i-\frac{1}{2}},x_{i+\frac{1}{2}})$ for $i\in \mathbb{Z}$ where $x_{i+\frac{1}{2}} = (i+\frac{1}{2})\Delta x$. Here, $\Delta x>0$ is a fixed mesh width. 
Moreover, we define the discrete time  points $t_n=n \Delta t$ for $n\in \mathbb{Z}^+$.  
As in any finite volume framework, we consider the cell averages of the solution $u$ given by  $u_j^n \approx
\frac{1}{\Delta x}\int_{I_j}u(t_n,x)dx$ at each 
$j\in \mathbb{Z}$.
We assume $\Delta t$ and $\Delta x$ fulfill the CFL condition 
$\Delta t /\Delta x \max_j|\sigma(A(u_j^n))|\leq C_{CFL}$ with $\sigma$ being the spectral radius.
Moreover, we assume that $x^n_\alpha$ is a numerical approximation to the $\alpha$-th jump position $x_{\alpha}(t_n)$. That position is obtained by using, e.g.,  shock capturing method or evaluation of the Rankine–Hugoniot condition, see e.g., \cite{HHNSS}. In the present work, we compute the jump position by discretizing the Rankine–Hugoniot condition.

Having this, we propose a simple way to compute the limits $u_{\alpha}^-$ and $u_{\alpha}^+$ at the discontinuity: (1) Choose the index $J_{\alpha}$ such that the jump position fulfills $x_{\alpha}^n\in I_{J_{\alpha}}$.
 (2) To reduce the effect of the numerical viscosity, we choose $K$ according to the mesh size. Then, we consider the interval  $[x_{\alpha}^n-K\Delta x, x_{\alpha}^n+K\Delta x]$. The values $u^\pm_{\alpha}$ are then defined by  $u^-_{\alpha}=u^n_{J_{\alpha}- K}$ and $u^+_{\alpha}=u^n_{J_{\alpha}+ K}$, respectively. We refer to \cite{HHNSS} for further details. {The 
 estimate of the width of the numerical viscosity layer $K$ is non--trivial and we are not aware of a theoretical result in this direction. In this work, we therefore compute the analytical solution $u$ and then select the parameter $K$ to match with the analytical solution. 
 }

\begin{algorithm}\label{algo1}
\caption{Computation of the GTV (single-shock case)}

\textbf{Given data:}\\
(1) Initial data for \eqref{original}: the cell average $u_j^0=\frac{1}{\Delta x}\int_{I_j}u(0,x)dx$ for $j\in \mathbb{Z}$.\\
(2) Initial data for the GTV: the initial shift $\xi^0=\xi(0)$ and the cell average $v^0_j=\frac{1}{\Delta x}\int_{I_j}v(0,x)dx$ for $j\in \mathbb{Z}$. 

From $u_j^0$, we reconstruct 
the derivatives $[u_x]_j^0$ and compute: $w^0_j = v^0_j + \xi^0 [u_x]_j^0$.

\textbf{Updating in time from $t_n$ to $t_{n+1}$:}

\textbf{Step 1:} Update $u_j^{n+1}$ by the conservative numerical scheme:
\begin{align*}
 u^{n+1}_j - u^n_j + \frac{\Delta t}{\Delta x} (F^n_{j+\frac{1}{2}} - F^n_{j-\frac{1}{2}} )=0.
\end{align*}
Here $F_{j+\frac{1}{2}}^n = F(u^n_{j},u^n_{j+1})$ and $F$ is the numerical flux function for \eqref{original}. \\[2mm]
\textbf{Step 2:}
Track the shock position $x_{\alpha}(t_n)$ and compute the limits $(u^-,u^+)$ and $(w^-,w^+)$ accordingly.
Then compute $\xi'$ from \eqref{eq-2.1} and update $\xi^{n+1}=\xi^n+\Delta t \xi'$.\\[2mm]
\textbf{Step 3:} Update $w^{n+1}_j$ by exploiting a conservative numerical scheme:
\begin{align*}
 w^{n+1}_j - w^n_j + \frac{\Delta t}{\Delta x} (\widehat{F}^n_{j+\frac{1}{2}} - \widehat{F}^n_{j-\frac{1}{2}} )=0.
\end{align*}
Here $\widehat{F}^n_{j+\frac{1}{2}} = \widehat{F}(w^n_{j},w^n_{j+1},u^n_{j},u^n_{j+1},\xi')$ and $\widehat{F}$ is the numerical flux function for \eqref{eq:W}. 
\textbf{Output:} 

We obtain $u_j^n$, $w_j^n$ and $\xi^n$ for each $j$ and $n$. At last, we reconstruct 
the derivatives $[u_x]_j^n$ and compute the variation $v^n_j = w^n_j - \xi^n [u_x]_j^n$.

\end{algorithm}

Some remarks are in order. 
\begin{remark} 
 We compute the tangent vector by exploiting a conservative numerical scheme for the equation \eqref{eq:W}. In this way, the interface condition \eqref{BCw} is preserved automatically and no explicit discretization for this condition is needed. 
\par 
    For simplicity, we only present the algorithm with a first--order explicit Euler time discretization and simple numerical flux $\widehat{F}_{j+\frac12}$. It may be possible to extend both.  
\par     
 Note that our main motivation is the computation of the tangent vector. Thus, we assume that the numerical solution for the original problem \eqref{original} can be computed sufficiently accurate using e.g. a finer mesh than for $v_j^n$. 

\end{remark}

\begin{remark}\label{Remark4.2}
    We would like to compare our method with the one developed in \cite{MR3881242,MR4102600} based on the delta-type source. For the case where $\partial_xu=0$ a.e., we know that $v=w$ a.e. Then the equation \eqref{eq:W} has the same form as the sensitivity equation with delta-type source term located at the discontinuity in \cite{MR3881242,MR4102600}. Indeed, the term $\xi'u_x$ in \eqref{eq:W} corresponds to a delta function. In this particular situation, the main difference lies in the numerical treatments of the source term, e.g., we directly utilize the property of conservative schemes and do not employ shock indicators in our method. For the case where $\partial_xu$ is not necessarily zero, unlike the method in \cite{MR3881242,MR4102600}, our strategy can be directly applied for any conservative schemes without further restriction. The last but not the least, our method also explicitly computes the jump position $x_{\alpha}(t)$ and its shift $\xi_{\alpha}$ by  discretizing the RH condition and \eqref{eq-2.1}, which are useful in the construction of the first-order variation \eqref{construction}.
\end{remark}

\subsection{The case of multiple discontinuities }\label{sec:4.2}
In the previous subsection, we propose a method to compute the generalized tangent vector for the case where only one jump occurs. 
In general, the solution $u$ has $N$ jumps at points $x_1(t)<x_2(t)<\cdots<x_N(t)$. We use $x=y+x_{\alpha}(t)$ and rewrite \eqref{11eq5} as: 
\begin{equation}\label{eq4.2}
\partial_t w_{\alpha} + \partial_x [A(u)w_{\alpha} - \xi_{\alpha}' u]=0,\qquad x\in D_{\alpha}.
\end{equation}
Here $D_{\alpha}\in \mathbb{R}$ is a neighborhood of $x_{\alpha}(t)$ which contains no point of  discontinuity of $u$ other than $x_{\alpha}(t)$.

We assume there is no interaction between different jumps. Thus one can choose $D_{\alpha}$ for $\alpha=1,2,...,N$ such that $\mathbb{R} = \cup_{\alpha}D_{\alpha}$. 
By directly using Algorithm \ref{algo1}, we compute \eqref{eq4.2} to get $w_{\alpha}$ in each $D_{\alpha}$. Next, we provide a simple method to glue all $w_{\alpha}$ together. 

In each time step $t_n$, we take $D_\alpha^n=\{j\in \mathbb{Z}~|~j_{\alpha}^{-}\leq j \leq j_{\alpha}^{+}\}$ for each $\alpha=1,2,...,N$. 
The indexes $j_{\alpha}^{-},~j_{\alpha}^{+}\in \mathbb{Z}$ and $j_{\alpha+1}^{-},j_{\alpha+1}^{+}\in \mathbb{Z}$ are chosen such that 
\begin{align*}
x_{j_{\alpha}^{-}-\frac{1}{2}}\leq x_{\alpha}(t_n) < x_{j_{\alpha}^{+}+\frac{1}{2}},\quad  
x_{j_{\alpha+1}^{-}-\frac{1}{2}}\leq x_{\alpha+1}(t_n) < x_{j_{\alpha+1}^{+}+\frac{1}{2}},\quad
j_{\alpha}^{+} - J = j_{\alpha+1}^{-} + J.
\end{align*}
where $J$ is a nonnegative integer. It is reasonable to choose such $D_{\alpha}^n$ at each time step $t^n$, since no interaction of shocks occurs (by assumption).
For the spatial grids $x_j\in D_{\alpha}^n$, we compute $w_{\alpha,j}^{n+1}$ by Algorithm \ref{algo1}. Notice that, at the point $j=j_{\alpha}^+$, the computation of the numerical flux depends on $w_{\alpha,j+1}^n$. Since $j+1\in D_{\alpha+1}^n$, we only have the value of $w_{\alpha+1,j+1}^n$. 
Motivated by the relation \eqref{w-variable}, we take
\begin{equation}\label{eq-glue}
w_{\alpha,j+1}^n = w_{\alpha+1,j+1}^n + (\xi^n_{\alpha} - \xi^n_{\alpha+1}) [u_x]^n_{j+1}
\end{equation}
at the point $j=j_{\alpha}^+$.

After the computation of $w_{\alpha}^n$ for $\alpha=1,2,...,N$, we take
$$
v^{n}_j = \left\{\begin{array}{ll}
    w_{\alpha,j}^{n}+ \xi^n_{\alpha} [u_x]^n_{j} &\quad j_{\alpha}^- + J \leq j < j_{\alpha}^+ -J \\[3mm]
    w_{\alpha+1,j}^{n}+ \xi^n_{\alpha+1} [u_x]^n_{j} &\quad j_{\alpha+1}^- + J \leq j < j_{\alpha+1}^+ -J. 
\end{array}\right.
$$

We end this section by presenting the numerical construction of the first-order variation. Due to the numerical viscosity, the numerical solution of $u$ will not exhibit a jump discontinuity. Thus, a direct computation of \eqref{construction} will lead to a spike near the jump point $x=x_{\alpha}(t_n)$. To avoid this phenomenon, we reconstruct the solution $u$ by the following procedure:
\begin{enumerate}
    \item  Take the index $J_{\alpha}$, such that the jump position $x_{\alpha}^n\in I_{J_{\alpha}}$ and take the limits $u^-_{\alpha}$, $u^+_{\alpha}$ by the method given in Subsection \ref{sec:4.1}.
    \item  Reconstruct the solution $\widetilde{u}^n_j$ by
\begin{equation}\label{construction-u}
\widetilde{u}^n_j =  u_j^n + \sum_{\alpha}(u_\alpha^--u_j^n) \chi_{[x_{\alpha}^n-K\Delta x,x_{\alpha}^n)} +  \sum_{\alpha} (u_\alpha^+-u_j^n) \chi_{[x_{\alpha}^n,x_{\alpha}^n+K\Delta x]}.
\end{equation}
Here,  $x \to \chi_{[a,b]}(x)$ is the characteristic function on the interval $[a,b].$ 
\end{enumerate}

Furthermore, we compute the first-order variation by 
\begin{equation}\label{construction-n}
[u_\epsilon]_j^n = \widetilde{u}_j^n +\epsilon v_j^n + \sum_{\xi^n_{\alpha}<0} (u_{\alpha}^+-u_{\alpha}^-)\chi_{[x_{\alpha}^n+\epsilon \xi_{\alpha}^n, x_{\alpha}^n]} - \sum_{\xi_{\alpha}^n>0} (u_{\alpha}^+-u_{\alpha}^-)\chi_{[x_{\alpha}^n,x_{\alpha}^n+\epsilon \xi_{\alpha}^n]}.
\end{equation}

\section{Numerical experiments}\label{sec:5}
In this section, we present numerical experiments to show the validity of the proposed algorithm. The scalar case and the systems case are considered for the example of  Burgers' equation and the $2\times 2$ system \eqref{psystem}, respectively.

In Algorithm \ref{algo1} we use the first-order  (global) Lax-Friedrichs fluxes 
\begin{align*}
F(u^n_{j},u^n_{j+1}) &= \frac{1}{2}\Big(f(u_{j}^n) + f(u_{j+1}^n) - L (u_{j+1}^n-u_j^n) \Big)\\
\widehat{F}(w^n_{j},w^n_{j+1},u^n_{j},u^n_{j+1},\xi') &= \frac{1}{2}\Big(A(u_{j}^n)w_j^n - \xi' u_j^n 
 + A(u_{j+1}^n)w_{j+1}^n - \xi' u_{j+1}^n - L (w_{j+1}^n-w_j^n) \Big).
\end{align*}
Here, $L=\max_j|\sigma(A(u_j^n))|$ and $\sigma$ is the spectral radius.

\subsection{Scalar case: Burgers' equation}\label{section5.1}
For Burger's equation 
$$
u_t + \left(\frac{u^2}{2}\right)_x = 0,
$$
we consider the example in \cite{BM} where the initial data is given by
$$
u(0,x) = (1+\epsilon)x\cdot \chi_{[0,1]}(x).
$$
The exact solution is 
\begin{equation}\label{Burgers:exa}
u^{\epsilon}(t,x) = \frac{(1+\epsilon)x}{1+(1+\epsilon)t}\cdot \chi_{[0,\sqrt{1+(1+\epsilon)t}]}(x).
\end{equation}
As derived in Section 4, the equations for $(w,\xi)$ read
$$
w_t + (uw -\xi u)_x = 0,\qquad     
\xi' = \frac{w_-+w_+}{2}.
$$
Here are some details about the computation. We choose the time step $\Delta t$ to satisfy the CFL condition with $C_{CFL}=0.1$. We partition the spatial interval $[0,2]$ into $N_x$ cells with $N_x=20000$. 
Moreover, in the implementation of Algorithm \ref{algo1}, we take the limits $(w^-,w^+)$, $(u^-,u^+)$ as discussed in Section 4.1 with a viscosity parameter of $K = 10$. We obtain the jump position $x_{\alpha}^n$ by computing the Rankine–Hugoniot condition
$$
x_{\alpha}^{n+1} = x_{\alpha}^n + \frac{\Delta t}{2 } (u^++u^-).
$$

We compute the results up to time $t=0.2$. 
We obtain the numerical results $x_{\alpha}(t)|_{t=0.2}=1.0953$ for the jump position and $\xi(t)|_{t=0.2}=0.0914$ for  the generalized tangent vector. 
To compare it with the theoretical result, we also compute $x_{\alpha}$ and $\xi$ from the exact solution \cite{BM}:
$$
x_{\alpha}^{exact} = \sqrt{1+t}=1.0954,\qquad 
\xi^{exact} = \frac{t}{2\sqrt{1+t}}=0.0913.
$$
From this computation, we see that the numerical solution for $\xi$ is very close to the exact value with a relative error
$e_{\xi} = |\xi-\xi^{exact}|/\xi^{exact} \approx 0.1\% $.

\begin{figure}[h]
\centering
\includegraphics[width=2.5in]{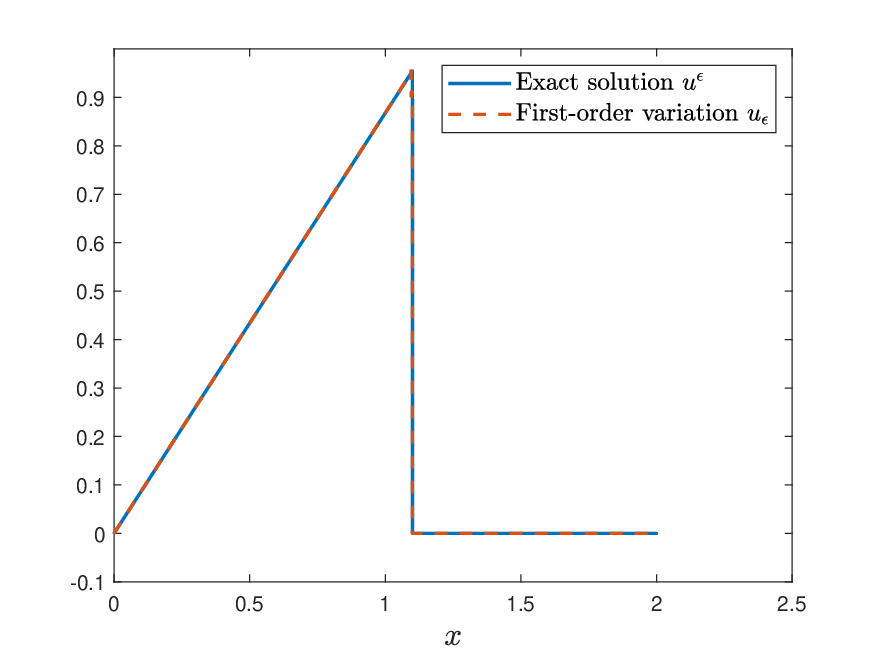}
\includegraphics[width=2.5in]{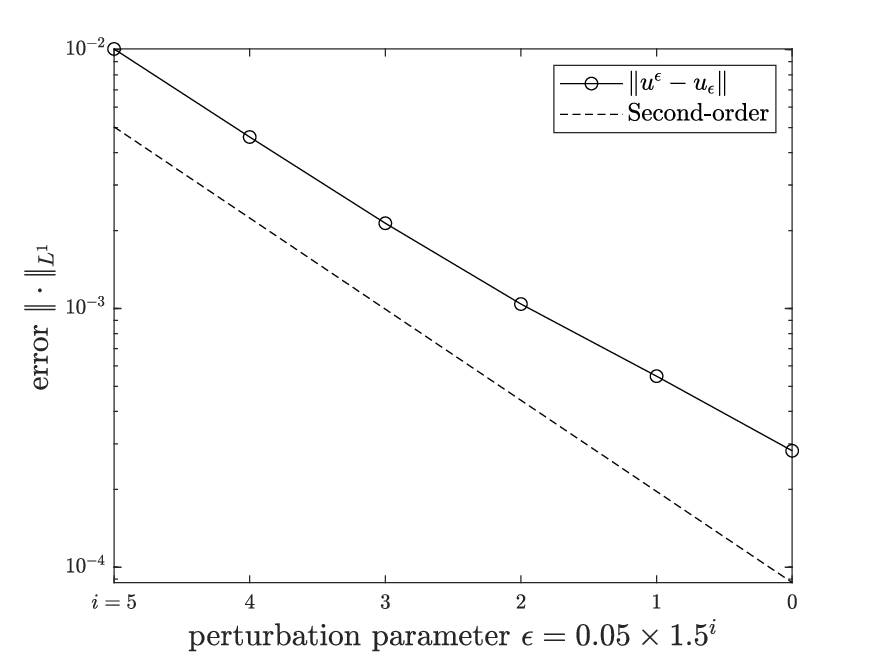}
\caption{Burger's equation: (left) the exact solution $u^{\epsilon}$ and the first-order variation $u_{\epsilon}$ computed from the GTV at $t=0.2$; (right) the difference $\|u^{\epsilon}-u_{\epsilon}\|_{L^1}$ between the exact solution and the first-order variation with different $\epsilon$.} 
\label{fig:E01}
\end{figure}

Having the numerical results for $u$, $w$, $x_{\alpha}$, $u^{\pm}$, and $\xi$ at $t=0.2$, we construct the first-order variation $u_{\epsilon}$ by
$$
u_{\epsilon}(x,t) = u(x,t) + \epsilon [w(x,t) - \xi u_x(x,t)] - [u^+-u^-] \chi_{[x_{\alpha}, x_{\alpha}+\epsilon \xi]}(x).
$$
The function $u$ is reconstructed according to \eqref{construction-u}. For the derivative $u_x$, we first compute $u_x^f$ and $u_x^b$ by $[u_x^f]_j^n=(u^n_{j+1}-u^n_j)/\Delta x$ and $[u_x^b]_j^n=(u^n_{j}-u^n_{j-1})/\Delta x$. Then we construct $u_x$ by 
$$
u_x = \left\{\begin{array}{lr}
   u_x^f  &  \quad |u_x^f|<|u_x^b|,\\[2mm]
   u_x^b  &  \quad |u_x^f|\geq|u_x^b|.
\end{array}\right.
$$
We plot the first-order variation and the exact solution with $\epsilon=0.05$ in Fig. \ref{fig:E01} (left). This experiment also validates that our method works for the case where $\partial_x u\neq 0$. The theoretical result in \cite{BM} indicates that the $L^1$-norm $\|(u^{\epsilon}-u_{\epsilon})/\epsilon\|_{L^1}$ should converge to $0$ as $\epsilon \rightarrow 0$. To verify this, we compute 
$\|u^{\epsilon}-u_{\epsilon}\|_{L^1}$ for $\epsilon = 0.05\times 1.5^i$ with $0\leq i\leq 5$ and plot the convergence result in Fig \ref{fig:E01} (right). 
{In this example,  we use the analytical solutions for $u^{\epsilon}$ in \eqref{Burgers:exa}.}
The dashed line represents the exact second-order convergence rate. From this figure, we see that the difference $\|u^{\epsilon}-u_{\epsilon}\|_{L^1}$ is of order $\epsilon^2$. This implies that $u_{\epsilon}$ is a valid first-order approximation for $u^{\epsilon}$ in the sense of $L^1$ norm.

\subsection{The $2\times 2$ system with single shock}\label{section5.2}
We consider the $2\times 2$ hyperbolic system in Section \ref{sec:3}:
\begin{eqnarray}\label{p-system}
    \left\{
    \begin{aligned}
        &\partial_t\rho + \partial_x q = 0,\\[1mm]
        &\partial_t q + \partial_x p(\rho) = 0,
    \end{aligned}\right.
\end{eqnarray}
where $p(\rho)=\kappa \rho^\gamma$. In what follows, we set $\kappa=1$ and $\gamma=2$. 
Consider the Riemann problem 
$$
(\rho^{\epsilon},~q^{\epsilon})(0,x) = \left\{\begin{array}{ll}
    (2+2\epsilon,~0), & x<0, \\[2mm]
    (1+\epsilon,~-\sqrt{3}(1+\epsilon)^{3/2}), & x\geq0.
\end{array}\right.
$$
A simple computation shows that the initial data is given on the curve \eqref{1-shock}. Therefore, for any $\epsilon$, only the 2-shock occurs and the exact solution is  
\begin{equation}\label{5.2exact}
(\rho^{\epsilon},~q^{\epsilon})(t,x) = \left\{\begin{array}{ll}
    (2+2\epsilon,~0), & x< \sqrt{3(1+\epsilon)} ~ t, \\[2mm]
    (1+\epsilon,~-\sqrt{3}(1+\epsilon)^{3/2}), & x\geq \sqrt{3(1+\epsilon)} ~t.
\end{array}\right.
\end{equation}

At first, we compute  the initial data for $(\xi,w)$. By the Taylor's expansion, we have $-\sqrt{3}(1+\epsilon)^{3/2}=-\sqrt{3} -\frac{3\sqrt{3}}{2}\epsilon+O(\epsilon^2)$. Then the initial data is 
$$
\xi(0) = 0,
\quad 
w(x,0) = \left\{\begin{array}{ll}
    (2,~0), & x<0, \\[2mm]
    \big(1,~-\frac{3\sqrt{3}}{2}\big), & x\geq0.
\end{array}\right.
$$
Now we use Algorithm \ref{algo1} to compute the generalized tangent vector $(\xi,v)$.
We take the CFL number as $C_{CFL}=0.5$ and partition the spatial interval $[-\pi,\pi]$ into $N_x$ cells with $N_x=20'000$. At the boundary,  the numerical flux is set to zero. 
Moreover, the jump position $x_{\alpha}^n$
is determined by solving the Rankine–Hugoniot condition
$$
x_{\alpha}^{n+1} = x_{\alpha}^n + \Delta t \frac{q^+-q^-}{\rho^+-\rho^-}.
$$
Here the limits $(\rho^{\pm},q^{\pm})$ and $(w_\rho^{\pm},w_q^{\pm})$ are evaluated following the procedure outlined in Section 4.1. In this context, the viscosity parameter is set to $K = 20$.

\begin{figure}[h]
\centering
\includegraphics[width=2.5in]{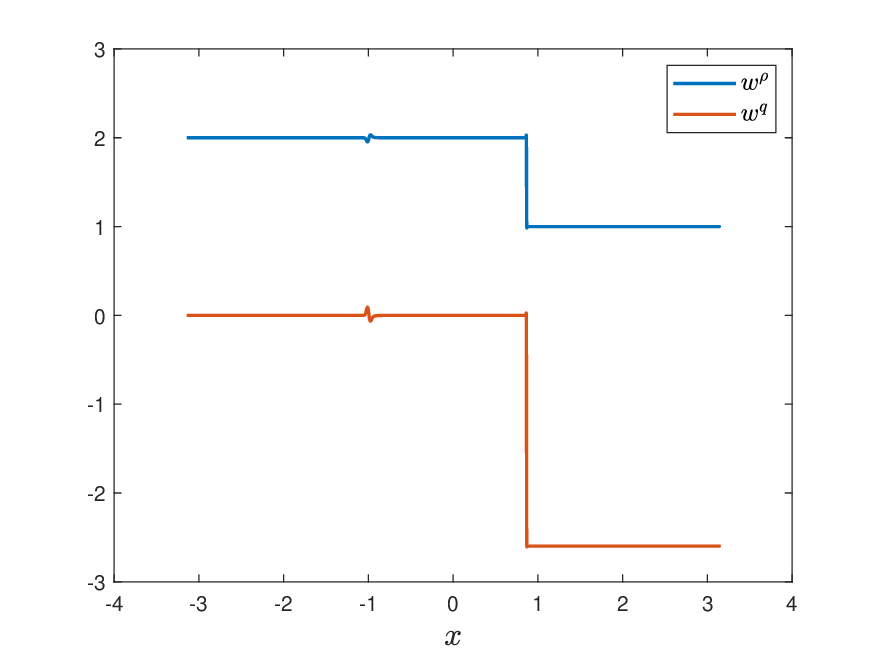}
\includegraphics[width=2.5in]{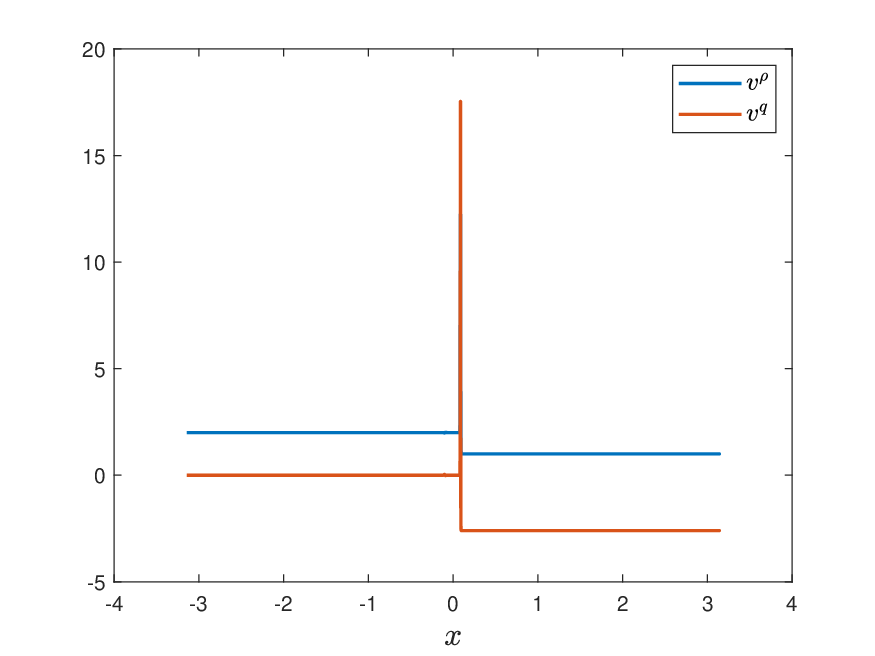}
\caption{Single shock case for the system \eqref{p-system}: (left) the numerical solution of $w=(w_\rho,w_q)$ at $t=0.5$ by using Algorithm \ref{algo1}. (right) the numerical solution of $v=(v_\rho,v_q)$ at $t=0.05$ by directly computing \eqref{theorem2.1:eq1}. This test shows the necessity to prescribe the correct interface condition.} 
\label{fig:E1-12}
\end{figure}

We compute the results at $t=0.5$. Namely, $x_{\alpha}|_{t=0.5}=0.8663$, $\xi|_{t=0.5}=0.4328$ and $w|_{t=0.5}=(w_\rho,w_q)|_{t=0.5}$ is shown in Fig. \ref{fig:E1-12} (left). From the exact solution \eqref{5.2exact}, we compute the theoretical results: 
$$
x_{\alpha}^{exact}|_{t=0.5} = 0.8660,\quad \xi^{exact}|_{t=0.5} = \lim_{\epsilon\rightarrow 0}\frac{\sqrt{3(1+\epsilon)}t-\sqrt{3}t}{\epsilon}\Big|_{t=0.5} = 0.4330.
$$
Thus we see that the numerical solution for $\xi$ is very close to the theoretical result with a relative error
$e_{\xi} = |\xi-\xi^{exact}|/\xi^{exact} \approx 0.05 \% $.

To show the importance of the interface condition \eqref{theorem2.1:eq3}, we also plot the numerical result by directly computing \eqref{theorem2.1:eq1} in Fig. \ref{fig:E1-12} (right). Namely, we use the first-order Lax-Friedrichs scheme to compute \eqref{theorem2.1:eq1} without considering the interface condition \eqref{theorem2.1:eq3}. At $t=0.05$, we observe the incorrect spike near the shock position. This test shows the necessity to prescribe the correct interface condition.

\begin{figure}[h]
\centering
\includegraphics[width=2.5in]{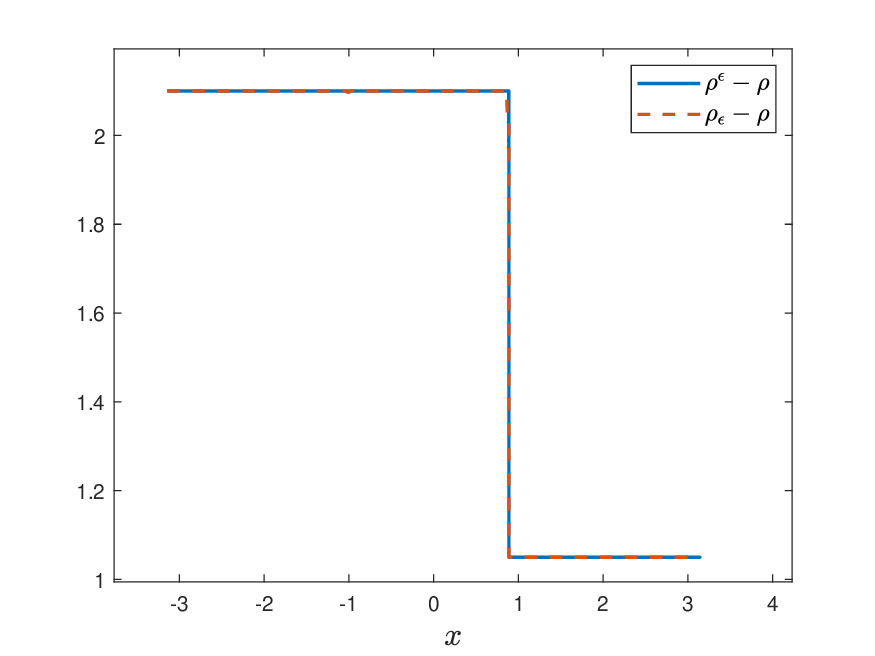}
\includegraphics[width=2.5in]{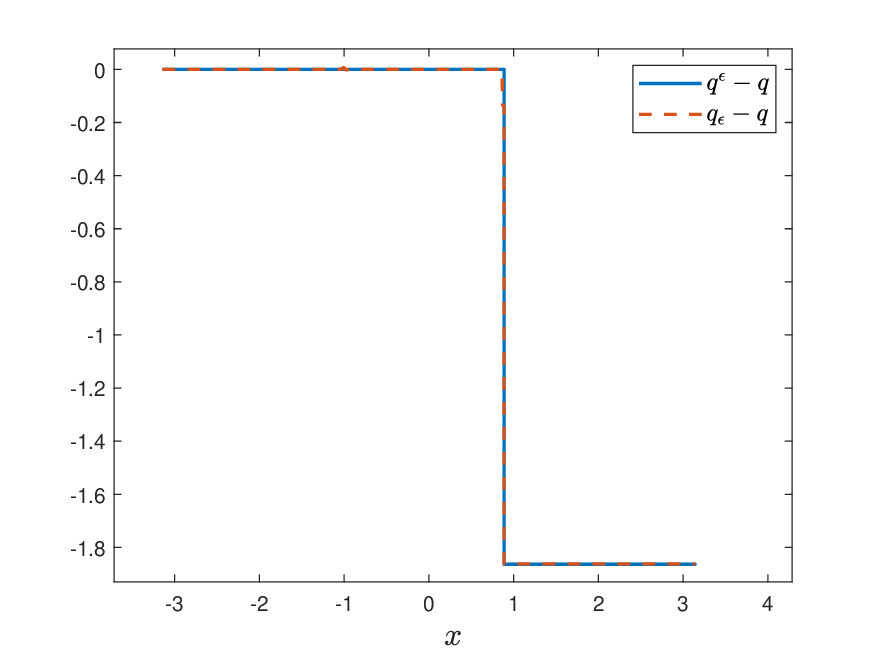}
\caption{Single shock case for the system \eqref{p-system}: the variation $(\rho^{\epsilon}-\rho^0,q^{\epsilon}-q^0)$ by solving the original problem and $(\rho_{\epsilon}-\rho^0,q_{\epsilon}-q^0)$ by computing the generalized tangent vector. } 
\label{fig2}
\end{figure}

\begin{figure}[h]
\centering
\includegraphics[width=2.5in]{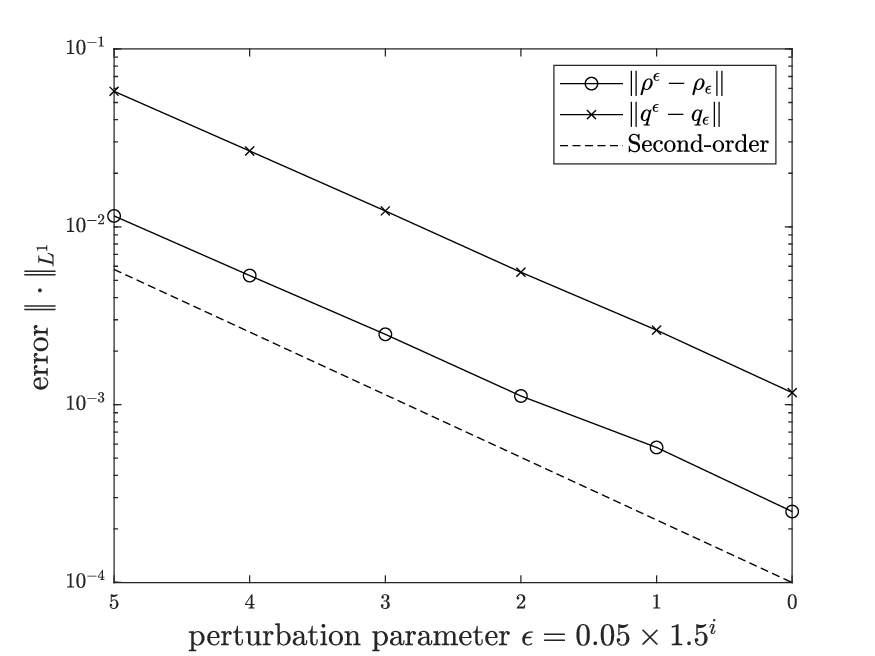}
\caption{Single shock case for the system \eqref{p-system}: the difference $\|(\rho^{\epsilon}-\rho_{\epsilon},q^{\epsilon}-q_{\epsilon})\|_{L^1}$ between the exact solution and the first-order variation.} 
\label{fig3}
\end{figure}

To verify that the computed generalized tangent vector is correct, we construct the first-order variation by 
\begin{equation}\label{5.2approx}
\begin{pmatrix}
\rho_\epsilon\\[1mm]
q_\epsilon   
\end{pmatrix}
= \begin{pmatrix}
\rho\\[1mm]
q 
\end{pmatrix}
+\epsilon 
\begin{pmatrix}
w_\rho\\[1mm]
w_q   
\end{pmatrix}
- \begin{pmatrix}
\rho^+-\rho^-\\[1mm]
q^+-q^- 
\end{pmatrix}
\chi_{[x_{\alpha}, x_{\alpha} + \epsilon \xi]}.
\end{equation}
Here $(\rho,q)$, $(w_{\rho},w_q)$, $(\rho^{\pm},q^{\pm})$, $x_{\alpha}$ and $\xi$ are given by the numerical results at $t=0.5$. We use the reconstruction of $(\rho,q)$ by \eqref{construction-u}.
Having this first-order variation, we compare it with the exact solution  $(\rho^\epsilon,q^\epsilon)(t,x)$ at $t=0.5$.
For $\epsilon = 0.05$, we plot in Fig. \ref{fig2} the first-order variation $(\rho_\epsilon,q_\epsilon)$ by \eqref{5.2approx} and the exact solution by \eqref{5.2exact}. 
Furthermore, we compute the $L^1$-norm of the difference $(\rho^{\epsilon}-\rho_{\epsilon},q^{\epsilon}-q_{\epsilon})$. Theoretically, we know that $\|(\rho^{\epsilon}-\rho_{\epsilon},q^{\epsilon}-q_{\epsilon})/\epsilon\|_{L^1}\rightarrow 0$ as $\epsilon \rightarrow 0$. With different $\epsilon$, we compute the convergence rate in Fig \ref{fig3}. 
{The  analytical solutions for $(\rho^{\epsilon},q^{\epsilon})$ are used in equation  \eqref{5.2exact}.}
Similar to the case of Burgers' equation,  the difference $\|(\rho^{\epsilon}-\rho_{\epsilon},q^{\epsilon}-q_{\epsilon})\|_{L^1}$ is of order $\epsilon^2$. This implies that $(\rho_{\epsilon},q_{\epsilon})$ is a  first-order approximation in  $L^1$.

\subsection{The $2\times 2$ system with two shocks}\label{section5.3}
In this subsection, we also consider the system \eqref{p-system} with $\kappa=1$ and $\gamma=2$. By properly choosing the initial data, we construct an example with two shocks. The initial data is given by
$$
(\rho^{0},q^{0})(0,x)=\left\{\begin{array}{ll}
     (2,~0) &\qquad x < x_1(0),\\[2mm]
     (2.241 ,-0.4963) &\qquad x_1(0) <x<x_2(0),\\[2mm]
     (1,-2.7304) & \qquad x > x_2(0).
\end{array}\right.
$$
Here $x_1(t)=-2.0593(t+0.1)$ and $x_2(t)=1.8002 (t+0.1)$.
Now we consider the perturbation of the initial data
$$
\begin{pmatrix}
\rho^{\epsilon}\\[1mm]
q^{\epsilon}
\end{pmatrix}(0,x) 
= 
\begin{pmatrix}
\rho^{0}\\[1mm]
q^{0}
\end{pmatrix}(0,x) + \epsilon \begin{pmatrix}
w_\rho\\[1mm]
w_q
\end{pmatrix}(0,x), \qquad x_i^{\epsilon}(0) = x_i(0) + \epsilon \xi_i(0),\quad i=1,2
$$
with $\xi_1(0)=-0.0728$, $\xi_2(0)=0.05554$ and
$$
(w_\rho,w_q)(0,x)=\left\{\begin{array}{ll}
     (2,-1.8839) &\qquad x < x_1(0),\\[2mm]
     (1,~0) &\qquad x_1(0) <x<x_2(0),\\[2mm]
     (1,-0.6893) & \qquad x > x_2(0).
\end{array}\right.
$$
With this initial condition, we compute the numerical solutions of $(\rho^{\epsilon},q^{\epsilon})$ for different $\epsilon$ by using the first-order finite volume scheme with Lax--Friedrichs flux. The spatial interval $[-\pi,\pi]$ is divided into $N_x$ cells with $N_x=20'000$.
On the other hand, we compute $w_\rho$ and $w_q$ 
by using Algorithm \ref{algo1} and the strategy in Section \ref{sec:4.2}. In this case, the value of $u_x$ in \eqref{eq-glue} is equivalent to zero.
This computation uses CFL number  $C_{CFL}=0.1$. The spatial interval $[-\pi,\pi]$ is also divided into $N_x$ cells with $N_x=20'000$. 
At the boundary, we take the numerical flux to be zero. We take the limits $(\rho^{\pm},q^{\pm})$ and $(w_\rho^{\pm},w_q^{\pm})$ at each time step by exploiting the method in Section \ref{sec:4.1} with the viscosity parameter $K=50$. 
Having $(w_\rho,w_q)$, we compute the first-order variation according to the formula \eqref{construction-n}.
Note that we also use the reconstruction of $u$ by \eqref{construction-u}.

\begin{figure}[h]
\centering
\includegraphics[width=2.5in]{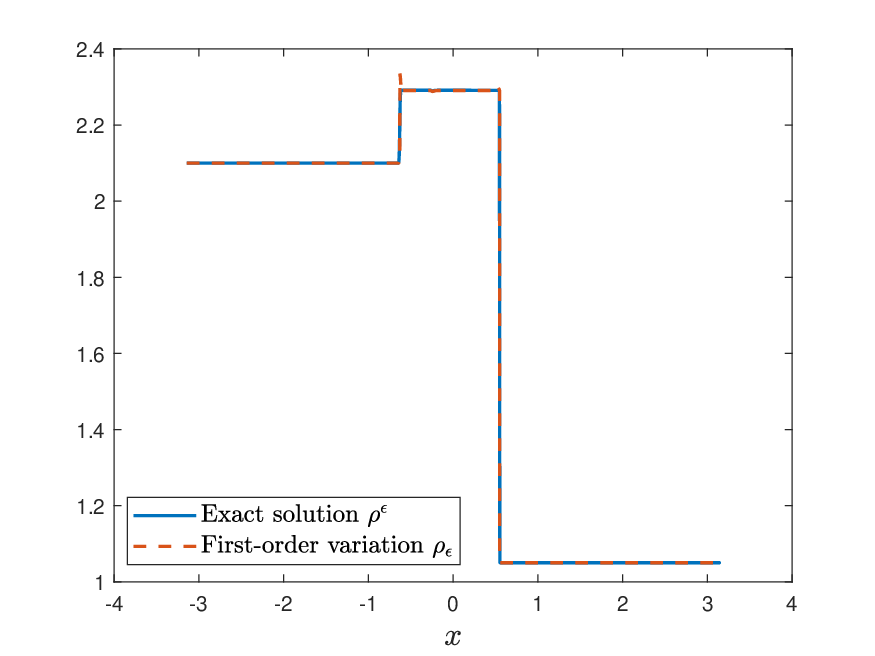}
\includegraphics[width=2.5in]{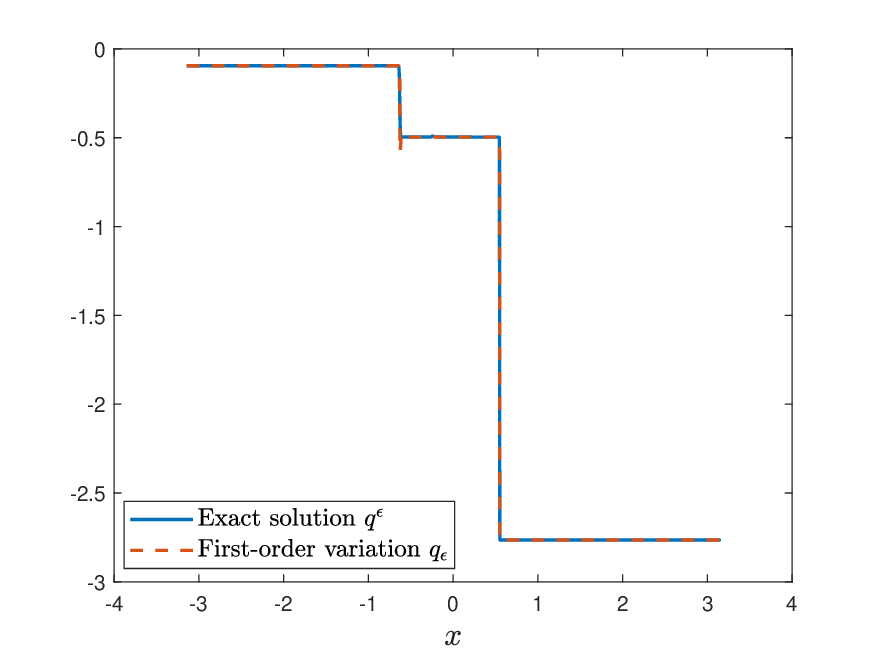}
\caption{
Two shocks case for the system \eqref{p-system}: The solution $(\rho^{\epsilon},q^{\epsilon})$ to the system \eqref{p-system} and the first-order variation $(\rho_{\epsilon},q_{\epsilon})$ computed from the GTV at $t=0.2$.
}
\label{psystem2}
\end{figure}

\begin{figure}[h]
\centering
\includegraphics[width=1.6in]{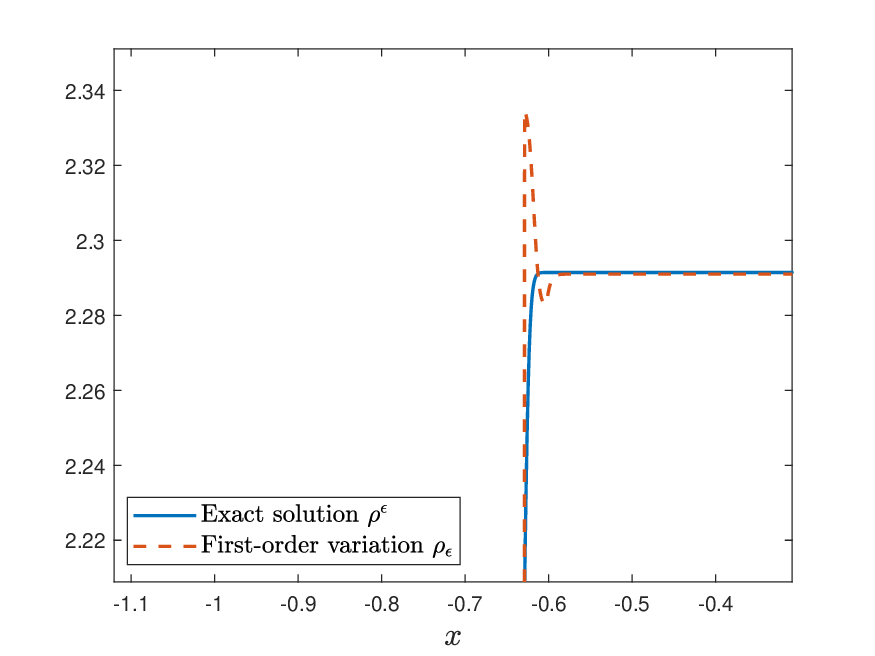}
\includegraphics[width=1.6in]{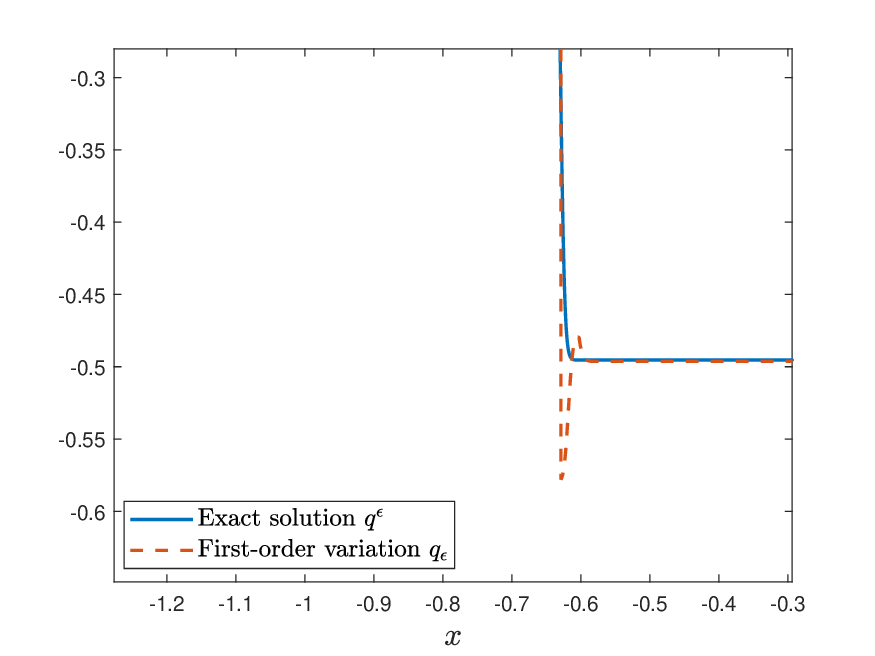}
\includegraphics[width=1.6in]{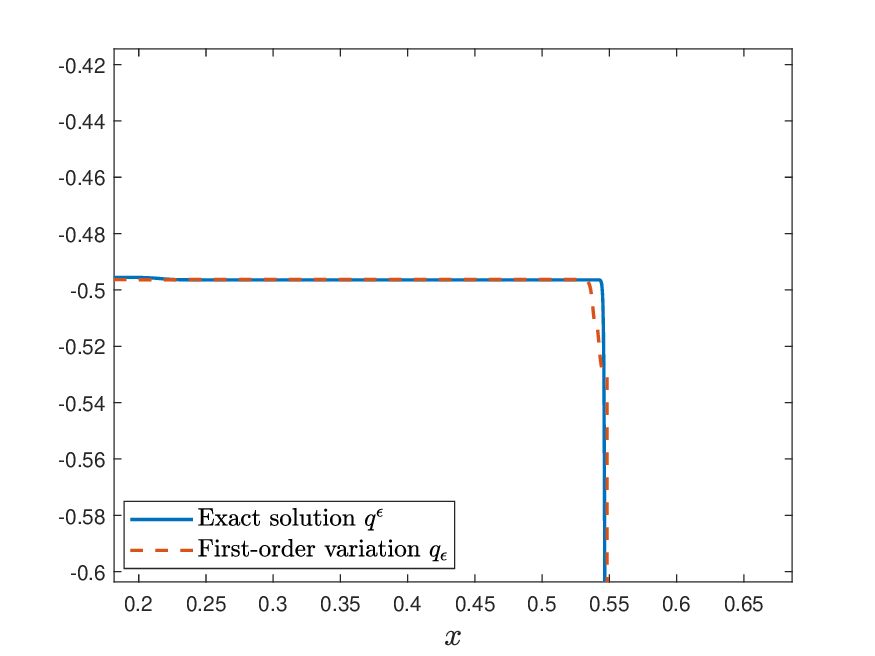}
\caption{Two shocks case for the system \eqref{p-system}: The solution $(\rho^{\epsilon},q^{\epsilon})$ to the system \eqref{p-system} and the first-order variation $(\rho_{\epsilon},q_{\epsilon})$ computed from the GTV at $t=0.2$. Details near $x=x_1(t)$ are shown in the first two figures and the detail near $x=x_2(t)$ is shown in the last figure.}
\label{fig6}
\end{figure}

In Fig \ref{psystem2}, we compare the computed first-order variation $(\rho_{\epsilon},q_{\epsilon})$ to the solution $(\rho^{\epsilon},q^{\epsilon})$ of the original system \eqref{p-system} with parameter $\epsilon=0.05$. 
Out of this figure, we observe that the first-order variation approximates the original solution well. 
There are small oscillations visible near the shock positions, see Fig. \ref{fig6} for zoomed-in details. In fact, this phenomenon also occurs for the theoretical result. More precisely, in the construction \eqref{construction}, the jump position for $u$ is shifted to $x_{\alpha}+\epsilon \xi_{\alpha}$ while the jump position for $v$ stays in $x_{\alpha}$.

For the perturbation of the jump position, the numerical result gives
\begin{equation}\label{Example3-xi1xi2}
\xi_1(t)|_{t=0.2}=-0.2184,\qquad \xi_2(t)|_{t=0.2}=0.1666.
\end{equation}
For this example, obtaining the exact results for $\xi_1$ and $\xi_2$ is not trivial. However, due to the  choice of the initial data, $(\rho^{\epsilon},q^{\epsilon})|_{t=0}$ are approximately given on the shock curves \eqref{1-shock} and \eqref{2-shock} with an error $O(\epsilon^2)$. Thus, we have the following approximate solution for sufficiently small $\epsilon$
$$
(\rho^{\epsilon},q^{\epsilon})(0,x)\approx\left\{\begin{array}{ll}
     (2+2\epsilon,-1.8839 \epsilon) &\qquad x < x_1^{\epsilon}(t),\\[2mm]
     (2.241+\epsilon,-0.4963) &\qquad x_1^{\epsilon}(t) <x<x_2^{\epsilon}(t),\\[2mm]
     (1+\epsilon,-2.7304-0.6893 \epsilon) & \qquad x > x_2^{\epsilon}(t),
\end{array}\right.
$$
where
$$
x_1^{\epsilon}(t) = \frac{-0.4963+1.8839\epsilon}{0.241-\epsilon} (t+0.1) = (-2.0593-0.728\epsilon)(t+0.1)+O(\epsilon^2),
$$
$$
x_2^{\epsilon}(t) = \frac{-2.2341-0.6893\epsilon}{-1.241} (t+0.1) = (1.8002+0.5554\epsilon)(t+0.1)+O(\epsilon^2).
$$
For comparison, we compute 
$$
\frac{x_1^{\epsilon}(t)-x_1(t)}{\epsilon}\Big|_{t=0.2} = -0.2184 + O(\epsilon),\qquad 
\frac{x_2^{\epsilon}(t)-x_2(t)}{\epsilon}\Big|_{t=0.2} =0.1666 + O(\epsilon).
$$
This matches the numerical result in \eqref{Example3-xi1xi2}.

\begin{figure}[h]
\centering
\includegraphics[width=2.5in]{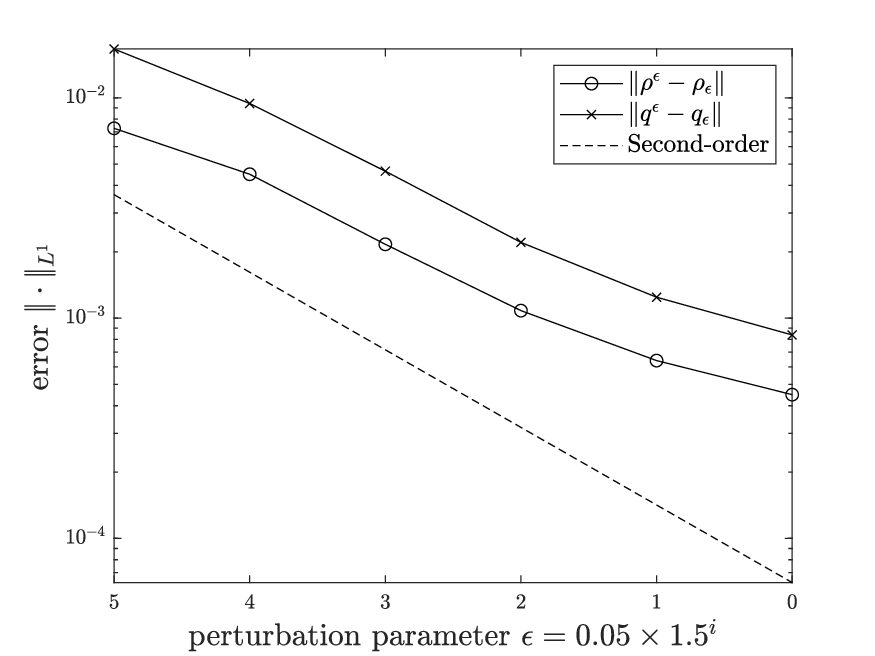}
\includegraphics[width=2.5in]{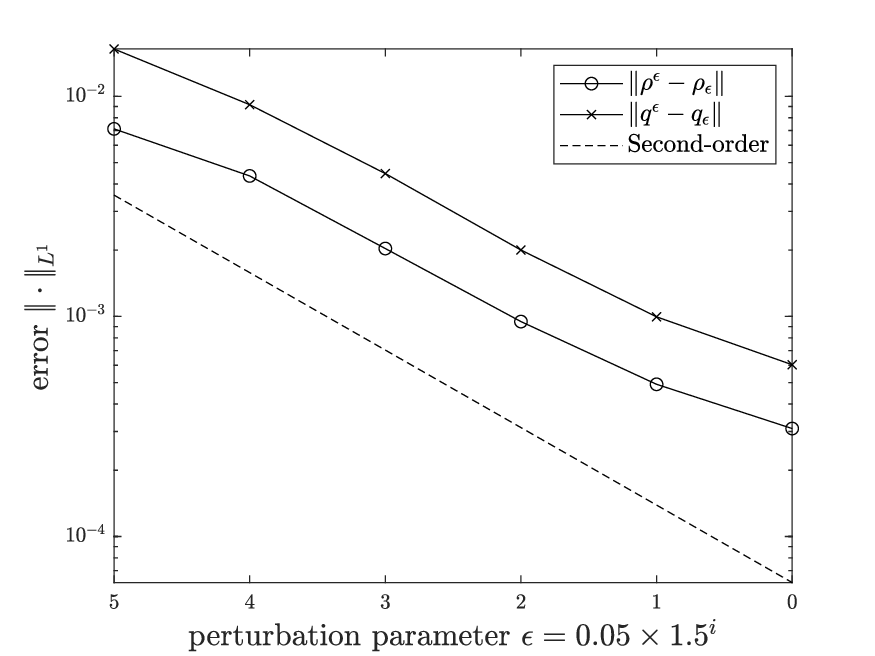}
\caption{Two shocks case for the system \eqref{p-system}: the difference $\|(\rho^{\epsilon}-\rho_{\epsilon},q^{\epsilon}-q_{\epsilon})\|_{L^1}$ between the exact solution and the first-order variation at $t=0.2$. Left: numerical solution with $N_x=20'000$, Right: numerical solution with $N_x=40'000$. } 
\label{psystem3}
\end{figure}

In order to validate the numerical approximation to the  tangent vector $(w_\rho,w_q)$ and $(\xi_1,\xi_2)$, we 
compute the convergence rate of $\|(\rho^{\epsilon}-\rho_{\epsilon},q^{\epsilon}-q_{\epsilon})\|_{L^1}$ in Fig \ref{psystem3} (left) with different $\epsilon$. Similarly to the previous cases, we see that the difference $\|(\rho^{\epsilon}-\rho_{\epsilon},q^{\epsilon}-q_{\epsilon})\|_{L^1}$ is approximately of order $\epsilon^2$ for sufficiently small $\epsilon$. This implies that $(\rho_{\epsilon},q_{\epsilon})$ is a  first-order approximation in  $L^1$. 
{
We observe that the convergence rate levels out as  $\epsilon$ becomes small. This is due to the numerical errors in solving the equations. When $\epsilon=0.05$, the error $\|(\rho^{\epsilon}-\rho_{\epsilon},q^{\epsilon}-q_{\epsilon})\|_{L^1}=O(10^{-4})$, which is of the order of the mesh size. To verify this, we also show the convergence result by using a finer grid with $N_x=40'000$.
}

\section{Outlook}

This work discusses a possible numerical computation of the generalized tangent vector and focuses on the numerical treatment for the interface conditions across discontinuities. We develop a simple numerical method in the framework of conservative schemes, enabling the automatic preservation of interface conditions without explicit discretization. A possible future work is to combine the strategy with other effective schemes, e.g. high-order Runge-Kutta discontinuous Galerkin scheme to reduce the numerical viscosity.

\begin{appendix}
\section{Details of the proof in Section 3}\label{AppendixA}

We give the proof of Lemma \ref{lemma12} by some tensor computations.
\begin{proof}
According to the definition of $l_i$, we have
\begin{equation}\label{lem32:eq1}
l_i(u^+,u^-)\bar{A}(u^+,u^-) = \lambda_i(u^+,u^-) l_i(u^+,u^-),\qquad i=1,2,...,n.
\end{equation}
The derivatives  $\partial/\partial u^+$ are computed on both sides.
Note that $\partial \bar{A}(u^+,u^-)/\partial u^+$ is a third-order tensor. To simplify the notation, we use the Einstein summation convention and write \eqref{lem32:eq1} as
$$
l_{ih}\bar{A}_{hj} = \lambda_i l_{ij},\qquad i,j=1,2,...,n.
$$
Taking derivative $\partial/\partial u^+_m$ on both sides ($u^+_m$ is the $m$-th component of $u^+$), we obtain
\begin{align*}
&\frac{\partial l_{ih}}{\partial u_m^{+}} \bar{A}_{hj} + l_{ih}\frac{\partial \bar{A}_{hj}}{\partial u_m^+}
=\frac{\partial \lambda_i}{\partial u_m^+} l_{ij} + \frac{\partial l_{ij}}{\partial u_m^+}\lambda_i.
\end{align*}
Multiplying the relation with $(u_j^+-u_j^-)$ and summing $j$ from $1$ to $n$, we have
\begin{align*}
&\frac{\partial l_{ih}}{\partial u_m^{+}} \bar{A}_{hj}(u_j^+-u_j^-) + l_{ih}\frac{\partial \bar{A}_{hj}}{\partial u_m^+}(u_j^+-u_j^-)
=\frac{\partial \lambda_i}{\partial u_m^+} l_{ij}(u_j^+-u_j^-) + \frac{\partial l_{ij}}{\partial u_m^+}\lambda_i (u_j^+-u_j^-).
\end{align*}
Since $\bar{A}_{hj}(u_j^+-u_j^-) = \lambda_{k_{\alpha}}(u_h^+-u_h^-)$, we have
\begin{align}\label{lemma32:eq2}
&\frac{\partial l_{ih}}{\partial u_m^{+}} (\lambda_{k_{\alpha}} - \lambda_i)(u_h^+-u_h^-) + l_{ih}\frac{\partial \bar{A}_{hj}}{\partial u_m^+}(u_j^+-u_j^-)
=\frac{\partial \lambda_i}{\partial u_m^+} (l_{ij}, u_j^+-u_j^-).
\end{align}

To simplify the second term, we compute
\begin{align*}
\frac{\partial \bar{A}_{hj}}{\partial u_m^{+}} 
= \int_0^1 \frac{\partial}{\partial u_m^{+}}\Big[A_{hj}(\theta u^+ + (1-\theta)u^-) \Big]d\theta
= \int_0^1 \theta \frac{\partial A_{hj}}{\partial u_m}(\theta u^+ + (1-\theta)u^-) d\theta. 
\end{align*}
Besides, by the definition of $A(u)$, we know that 
$$
\frac{\partial A_{hj}}{\partial u_m} = \frac{\partial^2 F_{h}}{\partial u_m \partial u_j} = \frac{\partial^2 F_{h}}{\partial u_j \partial u_m} = 
\frac{\partial A_{hm}}{\partial u_j}.
$$
Thus it follows that 
\begin{align*}
\frac{\partial \bar{A}_{hj}}{\partial u_m^{+}}(u_j^+-u_j^-) 
=&~ \int_0^1 \theta \frac{\partial A_{hm}}{\partial u_j}(\theta u^+ + (1-\theta)u^-)(u_j^+-u_j^-) d\theta\\[1mm]
=&~\int_0^1 \theta \frac{d}{d\theta}\Big[  A_{hm}(\theta u^+ + (1-\theta)u^-)\Big]\\[1mm]
=&~\theta A_{hm}(\theta u^+ + (1-\theta)u^-)\Big|_0^1 - \int_0^1  A_{hm}(\theta u^+ + (1-\theta)u^-)d\theta \\[1mm]
=&~A_{hm}(u^+) - \bar{A}_{hm}(u^+,u^-).
\end{align*}
Substituting this into \eqref{lemma32:eq2}, we have 
\begin{align*}
&\frac{\partial l_{ih}}{\partial u_m^{+}} (\lambda_{k_{\alpha}} - \lambda_i)(u_h^+-u_h^-) + l_{ih} [A_{hm}(u^+) - \bar{A}_{hm}(u^+,u^-)]
=\frac{\partial \lambda_i}{\partial u_m^+} (l_{ij}, u_j^+-u_j^-).
\end{align*}
This is the first relation in Lemma \ref{lemma12}. For the second case, the  proof is quite similar and we omit it.
\end{proof}

In the next, we clarify the notation $\langle D l_i(u^{+},u^{-})(w^{+},w^{-}), u^+-u^-\rangle$ in \eqref{BCw}:
\begin{align*}
Dl_{ih}(u^+,u^-)(w^{+},w^{-}) =&~ \lim_{\epsilon\rightarrow 0} \frac{l_{ih}(u^++\epsilon w^{+}, u^-+\epsilon w^{-})}{\epsilon} 
= \frac{\partial l_{ih}}{\partial u_m^+} w_m^+ + \frac{\partial l_{ih}}{\partial u_m^-} w_m^-.
\end{align*}
Thus we have
$$
\langle D l_i(u^{+},u^{-})(w^{+},w^{-}), u^+-u^-\rangle = \frac{\partial l_{ih}}{\partial u_m^+} w_m^+ (u^+_h-u^-_h) + \frac{\partial l_{ih}}{\partial u_m^-} w_m^- (u^+_h-u^-_h).
$$

\end{appendix}

\bibliographystyle{siamplain}

\end{document}